\documentclass[hidelinks,onefignum,onetabnum]{siamart251216}

\usepackage{amsfonts}
\usepackage{amssymb}
\usepackage{graphicx}

\newsiamremark{remark}{Remark}
\newsiamremark{assumption}{Assumption}
\AddToHook{env/remark/begin}{\crefalias{theorem}{remark}}
\AddToHook{env/assumption/begin}{\crefalias{theorem}{assumption}}
\crefname{assumption}{Assumption}{Assumptions}
\crefname{appendix}{Appendix}{Appendices}
\Crefname{appendix}{Appendix}{Appendices}

\usepackage{mathtools}
\usepackage{bm}
\usepackage{booktabs}
\usepackage[labelsep=period, textfont=it, labelfont={normalfont,sc}]{caption}
\usepackage{subcaption}
\usepackage{enumitem}

\newcommand{\eps}{\varepsilon}
\newcommand{\Reg}{\mathcal{R}}
\newcommand{\Jf}{\mathcal{J}}
\newcommand{\Ham}{\mathcal{H}}
\newcommand{\Hbar}{\overline{\mathcal{H}}}
\newcommand{\vhat}{\widehat{v}}
\newcommand{\grad}{\nabla}
\newcommand{\R}{\mathbb{R}}
\newcommand{\T}{\mathbb{T}}
\newcommand{\avg}[1]{\left\langle #1\right\rangle}
\newcommand{\Sph}{\mathbb{S}^1}
\newcommand{\norm}[1]{\left\lVert#1\right\rVert}
\newcommand{\abs}[1]{\left\lvert#1\right\rvert}
\newcommand{\argmin}{\operatorname*{arg\,min}}
\newcommand{\Sel}{\mathcal{D}}
\newcommand{\Imetric}{\mathcal{I}}

\headers{Decision Stability and Model Reduction for HJ Navigation}{R. Hu, B. Kong, and X. Yang}

\title{Decision Stability and Model Reduction for Hamilton--Jacobi Navigation in Multiscale Flows%
  \thanks{Submitted on \today.
    \funding{This work was partially supported by the ONR grant \#N00014-24-1-2432, the Simons Foundation (MP-TSM-00002783), and the NSF grant DMS-2420988.}}}

\author{Ruimeng Hu\thanks{Department of Mathematics and Department of
    Statistics and Applied Probability, University of California, Santa Barbara,
    CA 93106, USA (\email{rhu@ucsb.edu}).}
  \and
  Byungdoo Kong\thanks{Department of
    Statistics and Applied Probability, University of California, Santa Barbara,
    CA 93106, USA (\email{bkong@ucsb.edu}).}
  \and
  Xu Yang\thanks{Corresponding author. Department of Mathematics, University of California,
    Santa Barbara, CA 93106, USA (\email{xuyang@math.ucsb.edu}).}}

\usepackage{placeins}

\begin{document}
\maketitle

\begin{abstract}
Autonomous vehicles plan time-optimal routes through flow fields that are known only through reduced, finite-resolution models of an underlying multiscale flow. Classical Hamilton–Jacobi homogenization concerns convergence of the value function as the fine scale vanishes; but what a vehicle commits to is a route, and value accuracy need not imply route accuracy. We study this gap for a fixed-resolution reduction rather than in the homogenization limit. The committed route's dependence on the model-reduction error splits into two regimes. When the optimal route is a unique, nondegenerate minimizer with a regular value-function neighborhood, it depends Lipschitz-continuously on the reduction error, while the extra travel time from acting on the reduced model grows only quadratically with that error. When two routes are tied on the cut locus of the value function, the value remains stable, but an arbitrarily small change in the flow can switch the selected route, making route selection discontinuous in the reduction error. This dichotomy identifies where a reduced model must be accurate: not where the flow error is largest, but where the route choice is fragile. In the regular regime it yields a resolution-allocation principle, while an oracle marginal-regret benchmark illustrates the potential benefit of decision-aware resolution allocation relative to refinement based on flow-error magnitude.
\end{abstract}



\section{Introduction}\label{sec:intro}

An unmanned underwater or surface vehicle planning a minimum-time route through an ocean current does not know the true flow. It plans against a \emph{forecast}: a finite-resolution, data-assimilated model, such as an operational primitive-equation circulation product~\cite{lolla2014theory,lolla2014realistic}. The vehicle then moves through the true multiscale field, which the forecast only approximates. This discrepancy can change the committed route, so it cannot be judged by predicted arrival time alone.

A natural asymptotic framework for this problem is homogenization of Hamilton--Jacobi equations~\cite{lions1987homogenization,evans1992periodic}. The minimum time to the goal is the value function of a first-order Hamilton--Jacobi equation, and classical homogenization theory describes the effective value as the physical fine scale vanishes, with quantitative convergence rates available in convex settings~\cite{capuzzo2001rate,mitake2019rate}. The resulting effective problem is in general not obtained by simply averaging the flow. Classical homogenization therefore addresses a value-level asymptotic question. Here we instead ask a fixed-resolution question: when a vehicle plans with a reduced flow model, how faithfully does the resulting route represent the decision it would make in the true field?

Our central observation is that \emph{an accurate value does not imply an accurate decision}. The response of the committed route to model reduction splits into two regimes. In the regular regime, the optimal route is a unique, nondegenerate minimizer. There the route is stable under model perturbations: its displacement is first order in the model-reduction error, while the induced travel-time regret is second order. At the singular set associated with the cut locus~\cite{cannarsa2004semiconcave}, by contrast, multiple minimizing routes can compete. In the two-branch setting studied here, an arbitrarily small perturbation can reverse their ordering and switch the selected corridor by an $O(1)$ distance even while the value remains nearly unchanged. The general two-branch switching mechanism is analyzed in \Cref{subsec:singular}, and the opposing-jet experiment in \Cref{sec:sweep} illustrates it numerically. The accuracy required of a reduced model is therefore governed not simply by where the flow error is large, but by where the route decision is sensitive to that error.

This dichotomy suggests measuring model quality through its effect on the decision rather than through a global field-error norm. In the regular regime, the second-order route-regret sensitivity leads to an asymptotic principle for allocating model resolution. We also use an oracle marginal-regret benchmark to quantify the potential benefit of allocating resolution according to its effect on the navigation decision rather than according to flow-error magnitude alone. The oracle benchmark is a diagnostic rather than a practical adaptive algorithm; constructing and validating a computable surrogate remains open.

The analysis draws on several established strands of theory. Hamilton--Jacobi homogenization and its convergence rates~\cite{lions1987homogenization,evans1992periodic,capuzzo2001rate,mitake2019rate} provide the value-level multiscale background. The Zermelo--Randers formulation~\cite{bao2004zermelo} provides the geometric and variational description of time-optimal routes used here, while related Hamilton--Jacobi approaches address time-optimal navigation in dynamic flows~\cite{lolla2014theory} and reachability in differential games~\cite{mitchell2005time}. Our regular perturbation analysis is related to parametric stability of nondegenerate minimizers~\cite{dontchev2009implicit,bonnans2000perturbation}, and the allocation viewpoint is related in spirit to goal-oriented error estimation~\cite{becker2001optimal}, where approximation effort is directed toward a quantity of interest rather than a global field norm. The distinction here is that the object of interest is the selected route itself, whose dependence on the model becomes nonsmooth at route-selection boundaries. What is new is the decision-level analysis of multiscale reduction: the regular/singular stability dichotomy, quantitative route and regret sensitivity, and the resulting implications for where model resolution matters. 

The emphasis on committed routes also connects to our earlier work on inference and information control in navigation. In \cite{hu2026strategic}, observed route selection determines what can be learned about an agent's objective, while \cite{hu2026covertness} studies navigation in uncertain flows when the route itself reveals information to an adversarial observer. Related questions of information revelation, concealment, sequential inference, and deception arise more broadly in our work on dynamic and adversarial decision problems~\cite{hu2026stackelberg,hu2026information, zhou2026adversarial,zhou2025sunzi,kim2026deception,zhou2026optimal}.

\vspace{0.2em}

\noindent {\bf Main results.} The mathematical core is a dichotomy for the response of the committed route to model reduction. In the regular regime, \Cref{thm:regular} shows that the optimal route depends Lipschitz-continuously on the reduction error and that the resulting extra travel time is quadratic in that error, both for the committed route and for the executed feedback policy; the latter is established in \Cref{lem:track}. At a two-route tie, \Cref{prop:switch} shows that the value remains continuous while the selected route is discontinuous, so no single-valued route selection can be continuous through the tie. \Cref{cor:radius} quantifies the directional stability radius in the near-degenerate case, where two individually nondegenerate competing routes have a small positive action gap. These two regimes lead to distinct decision-sensitivity quantities. In the regular regime, \Cref{prop:alloc} converts the quadratic sensitivity into an asymptotic resolution-allocation principle. The numerical experiments illustrate both regimes and use an oracle marginal-regret benchmark to quantify the potential benefit of decision-aware refinement relative to refinement based on flow-error magnitude.

Section~\ref{sec:setup} sets up the model and its reduction, Section~\ref{sec:dichotomy} develops the regular/singular stability theory, and Section~\ref{sec:sensitivity} its implications for sensitivity and resolution allocation. Section~\ref{sec:numerics} presents numerical experiments, and Section~\ref{sec:disc} concludes.

\section{Multiscale Zermelo navigation and model reduction}
\label{sec:setup}

This section fixes the navigation model and the reduction setting. We introduce the dynamics and value function, distinguish the physical microscale from the computational resolution underlying the reduced model, recast the time-optimal problem in Randers--Finsler form to make the optimal route the object of analysis, and define the regret and route-selection quantities used throughout the paper.

\subsection{Time-optimal navigation and value function}
\label{subsec:dynamics}

Let $w:\Omega\to\R^2$ denote a steady flow field. On a bounded, convex domain $\Omega\subset\R^2$ with target set $\Gamma=\overline B_\rho(x_T)\subset\Omega$ (radius $\rho$, center $x_T$), fix an initial point $x_0\in\Omega\setminus\Gamma$. A vehicle with constant relative speed $s>0$ moves under a unit heading control $u\in\Sph$,
\begin{equation}\label{eq:dyn}
\dot X(t)=w(X(t))+s\,u(t),\qquad u(t)\in\Sph.
\end{equation}
A feedback policy $\pi:\Omega\to\Sph$ generates the trajectory
\[
\dot X^\pi=w(X^\pi)+s\,\pi(X^\pi),\qquad X^\pi(0)=x_0,
\]
with traversal time
\begin{equation}\label{eq:policy-cost}
J_w(x_0;\pi)=\inf\{t\ge 0:X^\pi(t)\in\Gamma\}.
\end{equation}

The corresponding minimum-time value is
\[
V_w(x_0)=\inf_\pi J_w(x_0;\pi).
\]
The state is constrained to $\overline\Omega$, reflecting the finite planning region; we take $\Omega$ large enough that the optimal trajectories considered below remain in its interior.

Under the strict subcriticality condition
\[
\sup_{x\in\Omega}\abs{w(x)}\le s-c_0<s,
\]
assumed for every flow field considered, the target is reachable, the Hamiltonian
\[
\Ham_w(x,p)=s\abs{p}-w(x)\cdot p
\]
is coercive and convex in $p$, and $V_w$ is the unique uniformly continuous constrained viscosity solution of the HJ equation~\cite{crandall1983viscosity,bardi1997optimal}, with the state-constraint condition on $\partial\Omega$,
\begin{equation}\label{eq:hj}
s\abs{\grad V_w(x)}-w(x)\cdot\grad V_w(x)=1\quad\text{in }\Omega\setminus\Gamma,\qquad V_w|_\Gamma=0.
\end{equation}
Where $V_w$ is differentiable, the optimal feedback is
\begin{equation}\label{eq:feedback}
\pi_w^*(x)=-\,\frac{\grad V_w(x)}{\abs{\grad V_w(x)}}.
\end{equation}

\subsection{True, homogenized, and reduced models}
\label{subsec:models}

We take the true flow to be genuinely multiscale,
\begin{equation}\label{eq:multiscale}
v^\eps(x)=v_0(x)+b\!\left(x,\tfrac{x}{\eps}\right),
\end{equation}
with $b(x,\cdot)$ mean-zero and periodic in the fast variable. Here $\eps$ is a \emph{physical} microscale, distinct from the computational resolution used to construct a reduced model. Three conceptually different objects must be distinguished:
\begin{enumerate}[leftmargin=*,label=(\roman*)]
\item the true field $v^\eps$ and its exact value $V_{v^\varepsilon}$;

\item the \emph{homogenized} Hamiltonian $\Hbar$ and value $V^0$, obtained from a cell problem. There need not exist an effective drift $v_{\rm eff}(x)$ such that
\[
\Hbar(x,p)=s\abs{p}-v_{\rm eff}(x)\cdot p\qquad\text{for all }p;
\]

\item a \emph{coarse reduction} $\vhat=\mathcal R_Hv^\eps$ at computational resolution $H$, such as a filtered field, subsampled field, or forecast product. Unlike the homogenized Hamiltonian in (ii), $\vhat$ is itself a drift field and corresponds to the specialization $w=\vhat$ of the planning problem in \Cref{subsec:dynamics}.
\end{enumerate}

Homogenization concerns the passage (i)$\to$(ii), whereas this paper studies the decision consequences of the practical reduction (i)$\to$(iii) at fixed $H$. The two reductions should not be conflated: the homogenized problem need not itself be of Zermelo form, while the practical reduced model remains within the drift-field class used for navigation. Separating the physical scale $\eps$ from the computational scale $H$ distinguishes model reduction from numerical discretization.

For the fixed-resolution analysis below, we suppress the superscript $\eps$ and write $v$ for the true field. Thus $v$ denotes the field in which the vehicle actually moves, whereas $\vhat$ denotes the reduced field available to the planner.

\subsection{Route-level variational formulation}
\label{subsec:randers}

The Hamilton--Jacobi formulation characterizes the minimum-time value, whereas our stability analysis concerns the optimal route itself. We therefore use the equivalent Zermelo--Randers variational formulation of the same time-optimal navigation problem.

For a strictly subcritical flow field $w$, consider admissible paths $\gamma:[0,1]\to\Omega$ satisfying $\gamma(0)=x_0$ and $\gamma(1)\in\Gamma$. The minimum-time value can be written as
\[
V_w(x_0)=\inf_\gamma \Jf_w(\gamma),
\]
where
\begin{equation}\label{eq:randers}
\Jf_w(\gamma)=\int_0^1 F_w(\gamma,\dot\gamma)\,\mathrm{d}r, \qquad F_w(x,y)=\frac{\sqrt{(w\cdot y)^2+(s^2-\abs{w}^2)\abs{y}^2}-w\cdot y}{s^2-\abs{w}^2}.
\end{equation}
This is the length functional of the Randers metric associated with Zermelo navigation~\cite{bao2004zermelo}. Strict subcriticality is precisely the condition under which $F_w$ defines a strongly convex Finsler metric. We denote the set of minimizing routes by
\[
\mathcal M_w=\argmin_\gamma\Jf_w(\gamma).
\]
Minimizing routes are identified up to orientation-preserving reparameterization. The length functional is invariant under reparameterization, which introduces an automatic tangential degeneracy into its second variation. For the perturbation analysis in \Cref{sec:dichotomy}, we therefore use the associated energy functional
\begin{equation}\label{eq:energy}
\mathcal E_w(\gamma)=\tfrac12\int_0^1F_w(\gamma,\dot\gamma)^2\,\mathrm{d}r.
\end{equation}
Every length minimizer admits a constant-speed parameterization that also minimizes $\mathcal E_w$ and is therefore a critical point of the energy; passing to $\mathcal E_w$ removes the reparameterization gauge. Since a minimizing route first reaches the target on $\partial\Gamma$, we use the endpoint-constrained variation space
\begin{equation}\label{eq:variation-space}
\mathcal H_\gamma=\{h\in H^1([0,1];\R^2):h(0)=0,\quad h(1)\in T_{\gamma(1)}\partial\Gamma\},
\end{equation}
with the $H^1$ norm; the terminal condition reflects that a minimizer meets $\partial\Gamma$ transversally. The energy index form on $\mathcal H_\gamma$, obtained by differentiating smooth variations in a local endpoint chart, includes the $\partial\Gamma$-curvature term; its bounded $H^1$ extension is the meaning of $\delta^2\mathcal E_w(\gamma)$ below. We do not claim the Finsler energy is twice Fr\'echet differentiable on the whole $H^1$ curve space. For the length index form we quotient out tangential reparameterizations $\mathcal T_\gamma=\{a\dot\gamma:a\in H^1_0(0,1)\}$, use the quotient norm $\norm{[h]}_{\perp}:=\inf_{k\in\mathcal T_\gamma}\norm{h+k}_{H^1}$, and write $\Imetric_w(h,h)$ for the resulting form. Taylor expansions use a $C^1$ chart of regular curves, with velocities bounded away from zero.

\subsection{Regret and route-selection error}
\label{subsec:regret}

The planner computes the reduced-model feedback $\pi_{\vhat}^*$ but executes that policy in the true field $v$. If the reduced model admits multiple optimal feedback policies, $\pi_{\vhat}^*$ denotes the particular policy selected by the planner. The resulting trajectory therefore satisfies
\[
\dot X(t)=v(X(t))+s\,\pi_{\vhat}^*(X(t)).
\]
The corresponding \emph{executed-policy regret} is
\begin{equation}\label{eq:regret}
\Reg_{\mathrm{pol}}(\vhat;v,x_0)=J_v(x_0;\pi_{\vhat}^*)-V_v(x_0)\ge 0,
\end{equation}
the extra traversal time incurred by using the reduced-model policy in the true field. We suppress $x_0$ when it is fixed. It measures additional travel time only.

For the variational analysis, we use the corresponding \emph{committed-route regret}
\begin{equation}\label{eq:regret-route}
\Reg_{\mathrm{route}}(\vhat;v)=\Jf_v(\gamma_{\vhat})-\min_\gamma\Jf_v(\gamma),\qquad \gamma_{\vhat}\in\mathcal M_{\vhat}.
\end{equation}
It measures the loss from committing to a route selected by the reduced model and evaluating that route in the true field. If the reduced model admits multiple optimal routes, $\gamma_{\vhat}$ denotes the particular route selected by the planner. The two regrets therefore play different roles: $\Reg_{\mathrm{route}}$ is the route-level variational quantity used in the analysis, whereas $\Reg_{\mathrm{pol}}$ is the operational quantity measured in the experiments. Under \Cref{as:regular} and for $\vhat$ sufficiently close to $v$, both satisfy
\[
\Reg_{\mathrm{route}}(\vhat;v)=O\!\left(\norm{\vhat-v}_{C^2}^2\right),\qquad \Reg_{\mathrm{pol}}(\vhat;v,x_0)=O\!\left(\norm{\vhat-v}_{C^2}^2\right)
\]
by \Cref{thm:regular,lem:track}.

Regret measures performance loss, but not whether the reduced model selects the same route as the true model. When the true minimizing route is unique, $\mathcal M_v=\{\gamma_v\}$, we therefore define the \emph{route-selection error}
\begin{equation}\label{eq:selerr}
\Sel(\vhat;v)=d(\gamma_{\vhat},\gamma_v),\qquad d(\gamma,\gamma')=\norm{\gamma-\gamma'}_{C^0([0,1])}.
\end{equation}
The distance is evaluated using the constant-speed parameterizations associated with \Cref{eq:energy}, so that the curves are compared in the same gauge. In the regular regime, \Cref{thm:regular} shows that $\Sel(\vhat;v)=O(\norm{\vhat-v}_{C^2})$: selection error is first order while regret is second order. At an exact tie, however, the true minimizing route is not unique and there is no canonical single-valued route-selection error relative to $v$; the resulting selection instability is analyzed in \Cref{sec:dichotomy}.

\begin{remark}\label{rem:losses}
The committed-route regret should be distinguished from the raw effect of changing the flow while holding a route fixed,
\[
\Jf_v(\gamma_{\vhat})-\Jf_{\vhat}(\gamma_{\vhat}).
\]
The latter generally contains a first-order contribution in $v-\vhat$. By contrast, $\Reg_{\mathrm{route}}$ compares two routes in the same true field $v$: since $\gamma_v$ is a critical point of $\Jf_v$, its first variation vanishes there, so together with the first-order response of the minimizing route the leading term is governed by the second variation, giving the quadratic scaling of \Cref{thm:regular}. The corresponding bound for the executed-policy regret needs additional regularity of the feedback field and is proved in \Cref{lem:track}.
\end{remark}

\section{Stability of optimal decisions}
\label{sec:dichotomy}

We now study how an optimal route responds to perturbations of the flow model. The answer depends on the geometry of the minimizing set. When the optimum is unique and nondegenerate, the minimizing route varies smoothly with the model and the resulting regret is second order. When two well-separated routes are tied, by contrast, the value remains stable while the globally selected route can switch discontinuously.

\begin{figure}[t]
\centering
\includegraphics[width=\textwidth]{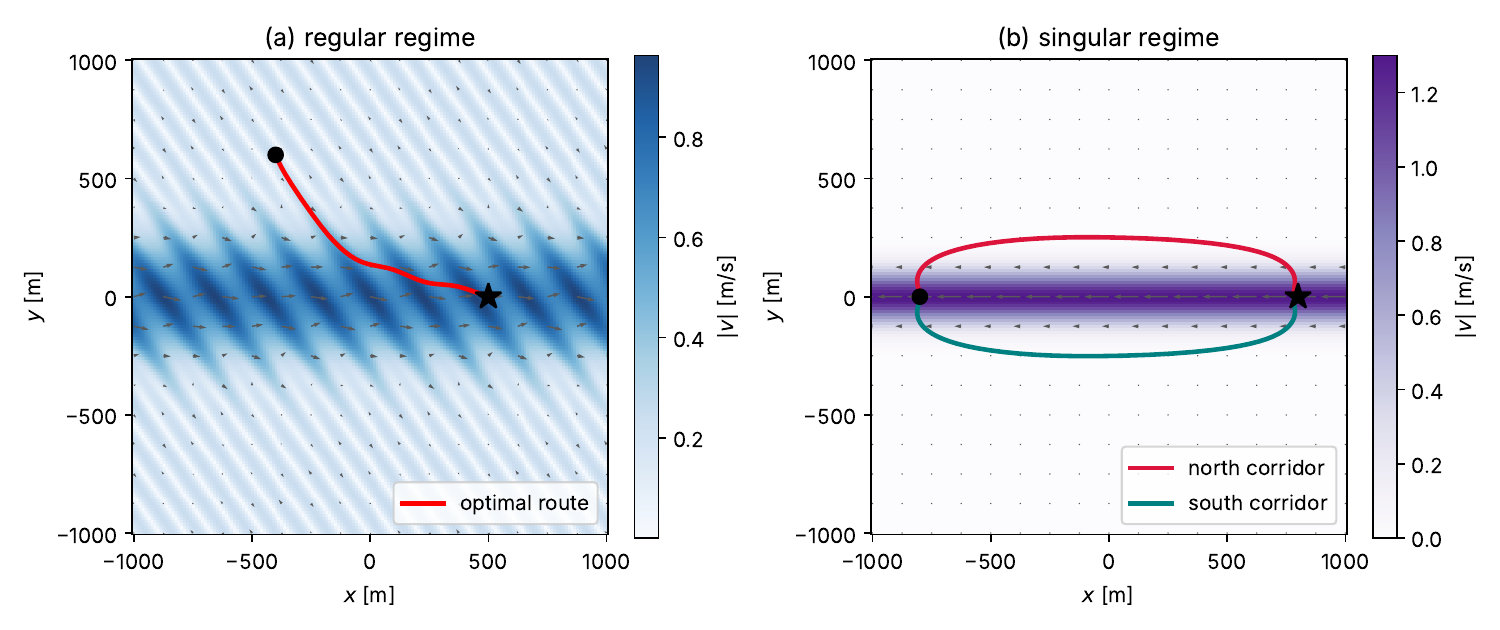}
\caption{Test flows. (a) Regular regime: shear jet with fine-scale oscillations; shading shows $|v|$, arrows the flow direction, and the red curve the optimal route from $x_0$ (dot) to the target (star). (b) Singular regime: symmetric opposing jet with the two competing corridors, computed with bias $\eta=\pm10\,$mm/s.}
\label{fig:setup}
\end{figure}

These two regimes are realized by our two test scenarios. In the regular example, a unique optimal route passes through a shear flow with fine-scale oscillations (its route is overlaid in \Cref{fig:setup}(a)). In the singular example, an opposing central jet creates competing northern and southern corridors (see \Cref{fig:setup}(b)); their minimizing routes are exactly tied in the symmetric case (\Cref{fig:switch}). We first clarify how this distinction appears in the homogenization limit and then develop the regular and singular perturbation theories. 

\subsection{Homogenization and route selection}
\label{subsec:homog}

Under the multiscale model \Cref{eq:multiscale} and the strict subcriticality condition, standard periodic homogenization results~\cite{lions1987homogenization,evans1992periodic} give
\[
V^\eps:=V_{v^\eps}\longrightarrow V^0
\]
locally uniformly, where $V^0$ solves
\[
\Hbar(x,\grad V^0)=1.
\]
Freezing the slow variable $x$ and writing $y=x/\eps\in\T^2$ for the fast variable, the effective Hamiltonian $\Hbar(x,p)$ is the constant for which the cell problem
\begin{equation}\label{eq:cell}
s\,\abs{p+\grad_y\varphi(y)}-\bigl(v_0(x)+b(x,y)\bigr)\cdot\bigl(p+\grad_y\varphi(y)\bigr)=\Hbar(x,p),\qquad y\in\T^2,
\end{equation}
admits a periodic viscosity corrector $\varphi=\varphi_{x,p}$.

The cell problem is not equivalent to simply averaging the drift. For example, take $v_0\equiv0$ and $b=c(y_1)e_1$ with $\avg{c}=0$. Here $e_1=(1,0)$, $y_1$ is the first component of the fast variable $y$, and $\avg{\cdot}$ denotes the average over one period. Then motion across one period in the $e_1$ direction has effective ground speed
\[
\left(\avg{(s+c)^{-1}}\right)^{-1},
\]
the harmonic rather than arithmetic mean. This is strictly below $s$ whenever $c\not\equiv0$, so even mean-zero microstructure can produce an $O(1)$ effect on effective travel time.

Value convergence alone does not determine the minimizing routes. Write $\Jf^\eps:=\Jf_{v^\eps}$. At the route level we assume that $\{\Jf^\eps\}$ is equicoercive and $\Gamma$-converges to an effective action $\Jf^0$ associated with $\Hbar$. Under these hypotheses, standard $\Gamma$-convergence theory~\cite{braides2002gamma,dalmaso1993gamma} implies that cluster points of minimizing routes belong to
\[
\mathcal M^0:=\argmin_\gamma\Jf^0.
\]
We use this minimizer-set convergence as a working hypothesis; establishing the corresponding $\Gamma$-limit for the Randers actions in \Cref{eq:randers} is outside our scope.

For each $\eps$, choose a minimizing route $\gamma^\eps\in\argmin_\gamma\Jf^\eps$. This gives a qualitative route-level dichotomy. If $\mathcal M^0=\{\gamma^0\}$, then the minimizing routes $\gamma^\eps$ converge to $\gamma^0$, so the homogenized problem identifies the limiting route as well as the limiting value. If $\#\mathcal M^0>1$, leading-order homogenization identifies only the minimizing set and does not determine which branch is selected by the finite-$\eps$ problem. Selection may then depend on higher-order microstructure not encoded by $\Hbar$.

\subsection{Regular regime: stable routes and quadratic regret}
\label{subsec:regular}

We first consider a unique nondegenerate minimizing route. The main result: route displacement is first order in the model perturbation, while both regrets are second order.

\begin{assumption}\label{as:regular}
For the start $x_0$ and target $\Gamma$, working in the constant-speed gauge of \Cref{eq:energy}, assume:
\begin{enumerate}[label=\textup{(A\arabic*)},leftmargin=3.0em]
\item $v\in C^2(\overline\Omega;\R^2)$ and
\[
\sup_\Omega\abs{v}\le s-c_0<s.
\]

\item The route $\gamma_v$ is the unique strict global minimizer of $\Jf_v$, its image is contained in the interior of $\Omega$, and it meets the target boundary $\partial\Gamma$ transversally at its terminal point. The value $V_v$ is $C^2$ in a neighborhood of $\gamma_v$ in $\overline\Omega\setminus\Gamma$, with one-sided $C^2$ regularity up to $\partial\Gamma$ at the terminal point; $C^2$ regularity \emph{across} $\partial\Gamma$ is not assumed and cannot be expected, as the distance function to $\Gamma$ already shows.

\item On $\mathcal H_{\gamma_v}$ of \Cref{eq:variation-space}, the energy index form is coercive: for some $c_E>0$,
\[
\delta^2\mathcal E_v(\gamma_v)[h,h]\ge c_E\norm{h}_{H^1}^2.
\]
The transverse length Hessian
\[
\Imetric_v:=\delta^2\Jf_v(\gamma_v)
\]
satisfies
\[
\Imetric_v(h,h)\ge\mu\norm{[h]}_{\perp}^2
\]
for some $\mu>0$. Geometrically, these assumptions exclude the conjugate or focal degeneracies relevant to the route perturbation problem.
\end{enumerate}
\end{assumption}

\begin{theorem}\label{thm:regular}
Under \Cref{as:regular}, there is a $C^2$-neighborhood $\mathcal U\ni v$ such that for every $\vhat\in\mathcal U$ the minimizing route $\gamma_{\vhat}$ is unique and depends $C^1$ on $\vhat$. Moreover,
\begin{equation}\label{eq:regular-bounds}
\begin{aligned}
\Sel(\vhat;v)&\le C\norm{\vhat-v}_{C^2},\\
0\le \Reg_{\mathrm{route}}(\vhat;v)&\le C\norm{\vhat-v}_{C^2}^2,\qquad 0\le \Reg_{\mathrm{pol}}(\vhat;v)\le C\norm{\vhat-v}_{C^2}^2.
\end{aligned}
\end{equation}
Along a one-parameter perturbation $\vhat=v+\eta\beta$, with $\beta\in C^2(\overline\Omega;\mathbb R^2)$ fixed and $|\eta|$ small enough that $\vhat\in\mathcal U$, the committed-route regret has the expansion 
\begin{equation}\label{eq:regular-quad}
\Reg_{\mathrm{route}}(v+\eta\beta;\,v)=\tfrac12\eta^2\Imetric_v(\xi_\beta,\xi_\beta)+o(\eta^2),\qquad\xi_\beta=D_v\gamma[\beta].
\end{equation}
Thus the route-selection error satisfies a linear bound in the model-reduction error, while both route and executed-policy regrets satisfy quadratic bounds.
\end{theorem}

The theorem also identifies the two local objects needed for the sensitivity theory of \Cref{sec:sensitivity}: the route-sensitivity operator $D_v\gamma$, which describes the first-order response of the route, and the transverse Hessian $\Imetric_v$, which converts that response into second-order regret. The proof uses a local persistence result that will also be needed branch by branch in the singular regime.

\begin{lemma}\label{lem:branch}
Let $v\in C^2(\overline\Omega;\R^2)$ be strictly subcritical and let $\gamma_0$ be a regular, constant-speed strict local minimizer of $\Jf_v$, contained in the interior of $\Omega$ and meeting $\partial\Gamma$ transversally at its first arrival. Suppose that its energy index form is coercive on $\mathcal H_{\gamma_0}$ and that its transverse length Hessian
\[
\Imetric_v:=\delta^2\Jf_v(\gamma_0)
\]
satisfies $\Imetric_v(h,h)\ge\mu\norm{[h]}_{\perp}^2$ for some $\mu>0$. Then there are a $C^2$-neighborhood $\mathcal U\ni v$ and a unique $C^1$ map $w\mapsto\gamma_w$ into $C^1([0,1];\R^2)$, with $\gamma_v=\gamma_0$, such that each $\gamma_w$ is a nondegenerate strict local minimizer of $\Jf_w$.

In a fixed local chart of curves whose terminal point belongs to $\partial\Gamma$, write
\[
\mathcal G(w,\gamma):=\delta\mathcal E_w(\gamma)
\]
for its first variation as an element of $\mathcal H_{\gamma_0}^*$. The bounded index-form operator $\mathsf J_v:\mathcal H_{\gamma_0}\to\mathcal H_{\gamma_0}^*$ is defined by $\langle\mathsf J_vh,k\rangle=\delta^2\mathcal E_v(\gamma_0)[h,k]$, including the endpoint-curvature contribution. The first-order sensitivity in the direction $\beta$ is
\begin{equation}\label{eq:lem-sensitivity}
\xi_\beta=D_v\gamma[\beta]=-\mathsf J_v^{-1}D_w\mathcal G(v,\gamma_0)[\beta]. \end{equation}
Moreover,
\[
d(\gamma_w,\gamma_0)\le C\norm{w-v}_{C^2},
\]
and
\begin{equation}\label{eq:lem-regret}
\Jf_v(\gamma_w)-\Jf_v(\gamma_0)=\tfrac12\Imetric_v\bigl(D_v\gamma[w-v],D_v\gamma[w-v]\bigr)+o\!\left(\norm{w-v}_{C^2}^2\right).
\end{equation}
No global minimality of $\gamma_0$ is required.
\end{lemma}

\begin{proof}
Use the energy Lagrangian $L_w=F_w^2/2$ and its Hamiltonian $K_w=(F_w^*)^2/2$, where $F_w^*(x,p)=s\abs{p}+w(x)\cdot p$. Strict subcriticality makes the Legendre map invertible away from zero. For $w\in C^2$, its Hamiltonian vector field is $C^1$ in phase space, and the flow depends $C^1$ on $w$ in the $C^2$ topology. We shoot on the fixed parameter interval $r\in[0,1]$. Let $(\gamma_{w,p_0},p_{w,p_0})$ solve Hamilton's equations with $\gamma_{w,p_0}(0)=x_0$ and $p_{w,p_0}(0)=p_0\in\R^2$. For the circle target set, put $b(z)=\abs{z-x_T}^2-\rho^2$ and choose a smooth nonvanishing tangent field $t_\Gamma(z)$ on an annular neighborhood of its boundary. Define the two-component residual
\begin{equation}\label{eq:shooting-residual}
S(w,p_0)=
\begin{pmatrix}
b(\gamma_{w,p_0}(1))\\
p_{w,p_0}(1)\cdot t_\Gamma(\gamma_{w,p_0}(1))
\end{pmatrix}\in\R^2.
\end{equation}
The first component fixes the terminal point on the target boundary; the second is the natural free-endpoint stationarity condition. Thus $S(v,p_0^\star)=0$ for the momentum generating $\gamma_0$. Geometric transverse arrival ensures that nearby solutions of these terminal conditions also have their first target arrival at $r=1$.

To prove that the $2\times2$ matrix $D_{p_0}S(v,p_0^\star)$ is invertible, take an initial-momentum increment in its kernel. Its Jacobi field $h$ has $h(0)=0$, has terminal value tangent to $\partial\Gamma$, and satisfies the linearized natural terminal condition, including differentiation of $t_\Gamma$. Integration by parts in the linearized Euler--Lagrange equation gives
\[
\delta^2\mathcal E_v(\gamma_0)[h,k]=0\qquad\text{for every }k\in\mathcal H_{\gamma_0}.
\]
Taking $k=h$ and using coercivity gives $h=0$. Invertibility of the Legendre map then makes the initial-momentum increment zero. Hence the matrix is invertible. The implicit-function theorem gives a unique $C^1$ map $p_0(w)$ near $p_0^\star$, and $\gamma_w=\gamma_{w,p_0(w)}$ is a $C^1$ branch in the constant-energy speed gauge. Local uniqueness here is among nearby stationary routes; it does not assert uniqueness of all local minimizers.

The index forms vary continuously along this $C^1$ branch, so energy coercivity and transverse length positivity persist. Taylor expansion in a $C^1$ endpoint chart gives strict local minimality. To identify the sensitivity, differentiate $\mathcal G(w,\gamma_w)=0$ in that same chart. This gives
\[
\langle\mathsf J_v\xi_\beta,k\rangle=-D_w\mathcal G(v,\gamma_0)[\beta](k)\qquad(k\in\mathcal H_{\gamma_0}).
\]
The bounded coercive index form makes $\mathsf J_v$ an isomorphism from $\mathcal H_{\gamma_0}$ to its dual, proving \Cref{eq:lem-sensitivity} without a $C^2$ assertion on the whole $H^1$ curve space. Flow dependence also gives $\norm{\gamma_w-\gamma_0}_{C^1}\le C\norm{w-v}_{C^2}$. The first variation of $\Jf_v$ vanishes at $\gamma_0$. Its second-order expansion in the regular $C^1$ chart, together with the $C^1$ branch expansion, proves \Cref{eq:lem-regret}. No global minimality of $\gamma_0$ is used.
\end{proof}

The next lemma connects the route-level perturbation to the policy actually executed in the true field.

\begin{lemma}\label{lem:track}
Let $X$ be the trajectory obtained by executing $\pi_{\vhat}^*$ in the true field $v$, with arrival time $T=J_v(x_0;\pi_{\vhat}^*)$. Suppose $T<\infty$ and $V_v$ is $C^1$ on an exterior neighborhood of $X([0,T])$, with one-sided regularity at its terminal point. Then the exact identity
\begin{equation}\label{eq:pol-identity}
1+\frac{\mathrm d}{\mathrm dt}V_v(X(t))=\frac{s}{2}\,\abs{\grad V_v(X(t))}\,\bigl|\pi_{\vhat}^*(X(t))-\pi_v^*(X(t))\bigr|^2
\end{equation}
holds almost everywhere on $0<t<T$, and consequently
\begin{equation}\label{eq:pol-regret}
\Reg_{\mathrm{pol}}(\vhat;v)=\frac{s}{2}\int_0^T\abs{\grad V_v(X(t))}\,\bigl|\pi_{\vhat}^*(X(t))-\pi_v^*(X(t))\bigr|^2\,\mathrm dt\;\ge\;0.
\end{equation}
Under \Cref{as:regular}, there is a fixed exterior tube $U$ around the reference route, including a target-boundary collar, on which
\begin{equation}\label{eq:feedback-stability}
\norm{\grad V_{\vhat}-\grad V_v}_{C^0(U)}+\norm{\pi_{\vhat}^*-\pi_v^*}_{C^0(U)}\le C\norm{\vhat-v}_{C^2}.
\end{equation}
For $\vhat$ sufficiently close to $v$, the executed trajectory stays in this regular neighborhood until arrival, and its arrival time is bounded uniformly. Thus the preceding identity applies and
\[
\Reg_{\mathrm{pol}}(\vhat;v)\le C\norm{\vhat-v}_{C^2}^2.
\]
\end{lemma}

\begin{proof}
Where $V_v$ is differentiable, $\pi_v^*=-\grad V_v/\abs{\grad V_v}$. Along the executed trajectory $\dot X=v(X)+s\,\pi_{\vhat}^*(X)$,
\[
\frac{\mathrm d}{\mathrm dt}V_v(X)=\grad V_v\cdot\bigl(v+s\,\pi_{\vhat}^*\bigr)=\grad V_v\cdot v+s\,\grad V_v\cdot\pi_{\vhat}^*.
\]
The stationary equation \Cref{eq:hj} gives $s\abs{\grad V_v}-v\cdot\grad V_v=1$, i.e. $1+\grad V_v\cdot v=s\abs{\grad V_v}$, while $\grad V_v=-\abs{\grad V_v}\,\pi_v^*$. Therefore
\[
1+\frac{\mathrm d}{\mathrm dt}V_v(X)=s\abs{\grad V_v}\bigl(1-\pi_v^*\cdot\pi_{\vhat}^*\bigr)=\frac{s}{2}\abs{\grad V_v}\,\bigl|\pi_{\vhat}^*-\pi_v^*\bigr|^2,
\]
using $\abs{\pi_{\vhat}^*-\pi_v^*}^2=2\,(1-\pi_v^*\cdot\pi_{\vhat}^*)$ for unit vectors; this is \Cref{eq:pol-identity}. Integrating from $0$ to $T$, with $X(0)=x_0$ and $V_v(X(T))=0$ on $\Gamma$, gives $T-V_v(x_0)=\Reg_{\mathrm{pol}}(\vhat;v)$ equal to the right-hand side of \Cref{eq:pol-regret}; nonnegativity is immediate.

For the feedback estimate \Cref{eq:feedback-stability}, persistence of a single route from \Cref{lem:branch} is not sufficient: one must construct the perturbed value $V_w$ and its feedback on a fixed exterior tube $U$, uniformly for $w$ near $v$, since (A2) controls the regularity of $V_v$ but not of every nearby $V_w$. This is carried out in \Cref{app:feedback}, which builds characteristic coordinates from the target boundary, proves their inverses exist on a fixed $U$ uniformly in $w$, and shows that the resulting characteristics are globally minimizing, so the local construction equals $V_w$. It yields $\norm{\grad V_{\vhat}-\grad V_v}_{C^0(U)}+\norm{\pi_{\vhat}^*-\pi_v^*}_{C^0(U)}\le C\norm{\vhat-v}_{C^2}$, which is \Cref{eq:feedback-stability}. The executed trajectory $\dot X=v(X)+s\,\pi_{\vhat}^*(X)$ is then a $C^0$-small perturbation of the true optimal trajectory $\dot X_v=v(X_v)+s\,\pi_v^*(X_v)$; since $v+s\,\pi_v^*$ is $C^1$ on $U$, Gronwall's inequality gives $\sup_t\abs{X-X_v}\le C\norm{\vhat-v}_{C^2}$ until arrival or lateral exit. The terminal normal velocity $v\cdot n-s\le-c_0<0$ (with $n$ the outward unit normal to $\partial\Gamma$ at the terminal point) drives $X_v$ strictly into $\Gamma$ through the collar, so for small perturbations $X$ arrives before leaving $U$, with arrival time bounded uniformly in $\vhat$. Substituting the $C^0$ feedback bound into \Cref{eq:pol-regret} yields the quadratic estimate; only $C^0$ control of the feedback difference is used.
\end{proof}

\begin{proof}[Proof of \Cref{thm:regular}]
\Cref{lem:branch}, applied with $\gamma_0=\gamma_v$, gives the unique local branch, the route-selection estimate, and the local quadratic regret expansion.

It remains to preserve global minimality. For a sequence $w_j\to v$ in $C^2$, let $\widetilde\gamma_j$ be global minimizers, parameterized with constant $F_{w_j}$ speed on $[0,1]$. Comparison with the local branch bounds their actions uniformly. Strict subcriticality gives
\[
F_w(x,y)\ge c\abs{y}\qquad(w\text{ near }v),
\]
and constant $F_{w_j}$ speed consequently gives a uniform Euclidean Lipschitz bound. A subsequence converges uniformly, and convexity in velocity gives lower semicontinuity of the limiting action. Moreover, on bounded-length competitors,
\[
\abs{\Jf_w(\gamma)-\Jf_v(\gamma)}\le C\norm{w-v}_{L^\infty}\operatorname{len}(\gamma).
\]
Here $\operatorname{len}(\gamma)$ denotes the Euclidean arclength of $\gamma$. Lower semicontinuity of the energy also identifies its constant-speed gauge: $\mathcal E_v$ at the limit is at most $\Jf_v(\gamma_v)^2/2$, while Cauchy--Schwarz gives the reverse inequality for every minimizing route. The limit is therefore a global minimizer of $\Jf_v$ in the same gauge, and uniqueness in (A2) identifies it as $\gamma_v$. It lies strictly inside $\Omega$, so the nearby minimizers also remain in the interior. The geodesic equation, together with the speed bound and uniform $C^2$ field bounds, then bounds their accelerations and gives convergence in $C^1$. Their initial momenta are consequently in the shooting neighborhood of \Cref{lem:branch}; its local uniqueness identifies $\widetilde\gamma_j$ with $\gamma_{w_j}$. By contradiction, every sufficiently nearby field has this unique global minimizing route. The committed-route estimates and \Cref{eq:regular-quad} follow from \Cref{eq:lem-regret}, and the executed-policy bound from \Cref{lem:track}.
\end{proof}

\subsection{Singular regime: discontinuous route selection}
\label{subsec:singular}

We next consider a two-route competition. At an exact tie, value stability no longer implies stability of the selected route: the value remains Lipschitz-continuous under perturbations, while the globally minimizing branch can switch discontinuously. A small positive action gap restores stability, but only up to a finite directional perturbation radius.

\begin{assumption}\label{as:singular}
Suppose
\[
\mathcal M_v=\{\gamma_1,\gamma_2\},\qquad d(\gamma_1,\gamma_2)\ge c>0.
\]
Each branch is nondegenerate in the sense of \Cref{lem:branch}, so it persists as a strict local minimizer under small perturbations. Moreover, there exist disjoint neighborhoods $N_i\ni\gamma_i$ with respect to $d$ and a constant $\kappa>0$ such that 
\[
\Jf_v(\gamma)\ge \Jf_v(\gamma_1)+\kappa
\]
for every admissible $\gamma\notin N_1\cup N_2$. Their actions are exactly tied,
\[
\Jf_v(\gamma_1)=\Jf_v(\gamma_2).
\]
The perturbation direction $\beta\in C^2(\overline\Omega;\mathbb R^2)$ separates their first-order branch values:
\begin{equation}\label{eq:branchfn}
G_\beta(\gamma_1)\neq G_\beta(\gamma_2),\qquad G_\beta(\gamma):=\left.\frac{\mathrm d}{\mathrm d\eta}\right|_{\eta=0}\Jf_{v+\eta\beta}(\gamma).
\end{equation}
\end{assumption}

We first record the stability of the value itself.

\begin{lemma}\label{lem:vstab}
Under strict subcriticality, there is $C>0$ such that for all admissible fields $v,w$,
\begin{equation}\label{eq:value-stability}
\norm{V_w-V_v}_{L^\infty(\Omega)}\le C\norm{w-v}_{L^\infty(\Omega)}.
\end{equation}
\end{lemma}

\begin{proof}
The Hamiltonians satisfy
\[
\abs{\Ham_w(x,p)-\Ham_v(x,p)}=\abs{(w-v)(x)\cdot p}\le\norm{w-v}_{L^\infty}\abs{p}.
\]
Strict subcriticality gives a uniform spatial Lipschitz bound for the corresponding minimum-time value functions. The standard comparison and continuous-dependence estimate for stationary Hamilton--Jacobi equations~\cite{crandall1992users,bardi1997optimal} then yield \Cref{eq:value-stability}.
\end{proof}

\begin{proposition}\label{prop:switch}
Under \Cref{as:singular}, let
\[
v_\eta=v+\eta\beta.
\]
For all sufficiently small $\eta$:
\begin{enumerate}[label=\textup{(\alph*)},leftmargin=2.2em]
\item for $\eta\neq0$, the global minimizer is unique and equals the branch with the smaller first-order perturbed action;

\item the value functions satisfy
\[
\norm{V_{v_\eta}-V_v}_{L^\infty(\Omega)}=O(\abs{\eta}),
\]
so the value converges uniformly through the branch switch;

\item after labeling the branches so that $G_\beta(\gamma_1)<G_\beta(\gamma_2)$,
\[
\gamma_{v_\eta}\longrightarrow\gamma_1\quad\text{as }\eta\to0^+,\qquad\gamma_{v_\eta}\longrightarrow\gamma_2\quad\text{as }\eta\to0^-.
\]
Consequently, no continuous single-valued selection of the global minimizer exists through $\eta=0$.
\end{enumerate}
\end{proposition}

\begin{proof}
By \Cref{lem:branch}, each $\gamma_i$ persists as a $C^1$ branch of strict local minimizers $\gamma_i^{(\eta)}$. Since the route variation contributes only at second order,
\[
\Jf_{v_\eta}(\gamma_i^{(\eta)})=\Jf_v(\gamma_i)+\eta G_\beta(\gamma_i)+O(\eta^2).
\]
Using the exact tie,
\[
\Jf_{v_\eta}(\gamma_1^{(\eta)})-\Jf_{v_\eta}(\gamma_2^{(\eta)})=\eta\bigl[G_\beta(\gamma_1)-G_\beta(\gamma_2)\bigr]+O(\eta^2).
\]
By \Cref{eq:branchfn}, this difference changes sign as $\eta$ crosses zero. The global action gap in \Cref{as:singular} persists under small perturbations by the same bounded-sublevel argument used in the proof of \Cref{thm:regular}, so no third route becomes globally minimizing. This proves (a). Part (b) follows from \Cref{lem:vstab}. Part (c) follows from the $C^1$ persistence of the two branches and $d(\gamma_1,\gamma_2)\ge c>0$.
\end{proof}

\begin{corollary}\label{cor:radius}
Label the branches so that $\gamma_1$ is the strict global minimizer, and let the branch-action gap be $\Delta(v)=\Jf_v(\gamma_2)-\Jf_v(\gamma_1)>0$. Suppose both branches are nondegenerate and, on a fixed interval $\abs{\eta}\le\eta_*$ on which they persist as in \Cref{lem:branch}, the perturbation $v_\eta=v+\eta\beta$ satisfies:
\begin{enumerate}[label=\textup{(\roman*)},leftmargin=2.2em]
\item a uniform slope separation $g:=G_\beta(\gamma_1)-G_\beta(\gamma_2)$ with $\abs{g}\ge g_0>0$, where $g_0$ is independent of $\Delta(v)$;
\item the branch-action difference $D(\eta):=\Jf_{v_\eta}(\gamma_1^{(\eta)})-\Jf_{v_\eta}(\gamma_2^{(\eta)})$ along the continued branches satisfies $D(\eta)=-\Delta(v)+\eta\,g+r(\eta)$ with a $C^1$ remainder obeying $\abs{r(\eta)}\le M\eta^2$ and $\abs{r'(\eta)}\le M\abs{\eta}$ on $\abs{\eta}\le\eta_*$, where $M\ge0$ and the interval is chosen so that $M\eta_*\le g_0/2$;
\item disjoint route neighborhoods $N_1,N_2$ with respect to $d$ in which the continued branches are the unique local minimizers, and a uniform gap outside their union: for some $\kappa>0$ and every $\abs{\eta}\le\eta_*$,
\[
\Jf_{v_\eta}(\gamma)\ge\min_{i=1,2}\Jf_{v_\eta}(\gamma_i^{(\eta)})+\kappa\quad\text{if }\gamma\notin N_1\cup N_2.
\]
\end{enumerate}
For a family with $\Delta(v)\to0$, require $g_0,M,\eta_*$ and $\kappa$ to be uniform across that family. Then there is $\Delta_0>0$ such that for $0<\Delta(v)<\Delta_0$ the global minimizer switches at a unique parameter $\eta_{\mathrm{crit}}$ with $\abs{\eta_{\mathrm{crit}}}\le\eta_*$, given by
\begin{equation}\label{eq:radius}
\eta_{\mathrm{crit}}=\frac{\Delta(v)}{G_\beta(\gamma_1)-G_\beta(\gamma_2)}+O\!\left(\Delta(v)^2\right),\qquad\abs{\eta_{\mathrm{crit}}}=\frac{\Delta(v)}{\abs{G_\beta(\gamma_1)-G_\beta(\gamma_2)}}+O\!\left(\Delta(v)^2\right).
\end{equation}
Thus the selected branch is stable along this perturbation direction up to a finite directional radius, which shrinks to zero as $\Delta(v)\to0$. The global minimizer is unique away from $\eta_{\mathrm{crit}}$, while the two continued branches are both minimizing there. The exact-tie discontinuity of \Cref{prop:switch} is recovered in this limit.
\end{corollary}

\begin{proof}
By \Cref{lem:branch} both branches persist as $C^1$ families $\gamma_i^{(\eta)}$ for $\abs{\eta}\le\eta_*$, and the branch-value expansion gives
\[
D(\eta)=\Jf_{v_\eta}(\gamma_1^{(\eta)})-\Jf_{v_\eta}(\gamma_2^{(\eta)})=-\Delta(v)+\eta\,g+r(\eta),
\]
with $g$ and $r$ as in (i)--(ii). On the entire stated interval,
\[
\operatorname{sign}(g)D'(\eta)\ge g_0-M\eta_*\ge g_0/2>0,
\]
so $D$ is strictly monotone with the sign of $g$. We have $D(0)=-\Delta(v)<0$, while at $\eta_a=2\Delta(v)/g$,
\[
D\!\left(\tfrac{2\Delta(v)}{g}\right)=\Delta(v)+r\!\left(\tfrac{2\Delta(v)}{g}\right)\ge \Delta(v)-\frac{4M}{g_0^2}\Delta(v)^2>0
\]
for $0<\Delta(v)<\Delta_0:=\min\{g_0^2/(8M),\,g_0\eta_*/2\}$; the first entry is interpreted as $+\infty$ when $M=0$. This also ensures $\abs{\eta_a}\le\eta_*$. Hence there is exactly one zero between $0$ and $\eta_a$, for either sign of $g$, and no other zero on $\abs{\eta}\le\eta_*$. At that zero,
\[
\left|\eta_{\mathrm{crit}}-\frac{\Delta(v)}{g}\right|=\frac{\abs{r(\eta_{\mathrm{crit}})}}{\abs{g}}\le\frac{4M}{g_0^3}\Delta(v)^2,
\]
which proves \Cref{eq:radius} with a uniform remainder. Condition (iii) excludes global minimizers outside $N_1\cup N_2$; within each neighborhood, local uniqueness identifies any global minimizer with its continued branch. Thus the change in the sign of $D$ is precisely the global branch switch.
\end{proof}

The regular and singular results thus identify different local mechanisms of decision sensitivity, which form the basis of the localized sensitivity measures introduced in \Cref{sec:sensitivity}. Numerical illustrations of regular regret scaling and of exact and near-degenerate branch switching are presented in \Cref{sec:numerics}.

\section{Decision sensitivity and resolution allocation}
\label{sec:sensitivity}

The stability results of \Cref{sec:dichotomy} suggest measuring reduction quality by its effect on the navigation decision rather than by field error alone. They also identify two different mechanisms of decision sensitivity. In the regular regime, the local response is governed by the route derivative $D_v\gamma$ and the transverse length Hessian $\Imetric_v$. Near a competition between distinct routes, stability is instead controlled by the branch-action gap and by how perturbations change the relative branch values. We use these objects first to define localized sensitivity measures and then to derive a regular-regime principle for allocating model resolution.

\subsection{Decision adequacy}
\label{subsec:adequacy}

We distinguish small performance loss from stability of the selected route.

\begin{definition}
For tolerances $\delta,\zeta\ge0$, a reduction $\vhat$ is \emph{$\delta$-route-regret-adequate} for a mission if
\[
\Reg_{\mathrm{route}}(\vhat;v)\le\delta.
\]
When the true minimizing route is unique, $\mathcal M_v=\{\gamma_v\}$, we call $\vhat$ \emph{$(\delta, \zeta)$-decision-adequate} if, in addition, $\Sel(\vhat;v)\le \zeta$.
\end{definition}

Decision adequacy captures information that route regret alone does not. Consider two well-separated, nondegenerate branches with a small true-field action gap $\Delta(v)=\Jf_v(\gamma_2)-\Jf_v(\gamma_1)>0$, where $\gamma_1$ is the global minimizer. By \Cref{cor:radius}, a perturbation can reverse the branch ranking once the directional stability radius is crossed, changing the selected route by $O(1)$ even though the route regret stays small. At the limiting exact tie both branches are optimal and no canonical single-valued route-selection error exists, as described in \Cref{prop:switch}.

This leads to an inverse question---given a tolerance on the navigation decision, how much model resolution is required, and where?

\subsection{Localized decision sensitivity}
\label{subsec:localized-sensitivity}

\paragraph{Regular regime}
In the regular regime, \Cref{thm:regular} identifies the route derivative $D_v\gamma$ and the transverse length Hessian $\Imetric_v$. Together they define the symmetric bilinear form
\begin{equation}\label{eq:Bform}
B_v(b_1,b_2)=\Imetric_v\bigl(D_v\gamma[b_1],D_v\gamma[b_2]\bigr).
\end{equation}
By \Cref{eq:regular-quad},
\[
\Reg_{\mathrm{route}}(v+\eta b;\,v)=\tfrac12\eta^2 B_v(b,b)+o(\eta^2).
\]

For a prescribed family of localized perturbation directions $\beta_K$, supported near a cell $K$ and normalized by $\norm{\beta_K}_{C^2}=1$, we define the \emph{probe-based local quadratic sensitivity}
\begin{equation}\label{eq:Q}
Q(K)=\lim_{\eta\to0}\frac{\Reg_{\mathrm{route}}(v+\eta\beta_K;\,v)}{\eta^2}=\tfrac12 B_v(\beta_K,\beta_K).
\end{equation}
Thus $Q(K)$ measures the second-order regret response to a localized probe; it depends on the probe family and is not an intrinsic spatial density. The off-diagonal $B_v(\beta_K,\beta_{K'})$ measures interaction between cells, which reappears in the allocation problem.

\paragraph{Near-singular regime}
Near a competition between two nondegenerate branches, the branchwise Jacobi response may remain regular even while the global route selection becomes fragile. The relevant quantities are instead the branch-action gap $\Delta(v)$ introduced in \Cref{cor:radius} and the difference between the first-order branch responses. For the localized direction $\beta_K$, define
\begin{equation}\label{eq:Sbranch}
S_{\mathrm{branch}}(K)=\frac{\abs{G_{\beta_K}(\gamma_1)-G_{\beta_K}(\gamma_2)}}{\Delta(v)}.
\end{equation}
The denominator measures proximity to a global branch tie, while the numerator measures how strongly perturbing cell $K$ changes the branch ranking. For normalized $\beta_K$ with nonzero slope separation and satisfying the uniformity and interval conditions of \Cref{cor:radius}, \Cref{eq:radius} identifies $S_{\mathrm{branch}}(K)$ as the reciprocal of the leading-order directional stability radius. If the numerator vanishes, the indicator detects no first-order separation; it does not exclude switching caused by higher-order terms. At an exact tie, $\Delta(v)=0$, this normalized indicator is undefined, reflecting the discontinuous selection described in \Cref{prop:switch}. Thus $Q(K)$ measures the quadratic response of a single persistent branch, whereas $S_{\mathrm{branch}}(K)$ measures the fragility of the ranking between competing branches.

\subsection{Resolution allocation in the regular regime}
\label{subsec:allocation}

We next ask how the regular-regime sensitivity can guide the allocation of model resolution. Let $e=\vhat-v$ denote the reduction error. Fix a smooth partition of unity $\{\chi_K\}$ subordinate to slightly enlarged cells and write
\[
e=\sum_K e_K,\qquad e_K=\chi_K e.
\]
For $\norm{e}_{C^2}$ small, \Cref{thm:regular,eq:Bform} give
\begin{equation}\label{eq:regret-full}
\Reg_{\mathrm{route}}(\vhat;v)=\tfrac12 B_v(e,e)+o\!\left(\norm{e}_{C^2}^2\right)=\tfrac12\sum_{K,L}B_v(e_K,e_L)+o\!\left(\norm{e}_{C^2}^2\right).
\end{equation}
The quadratic form is not diagonal in the cell decomposition: cross-cell terms $B_v(e_K,e_L)$ with $K\neq L$ need not vanish or have a fixed sign. Consequently, resolving two cells need not reduce regret by the sum of their individual effects.

To obtain a tractable allocation principle, define the cellwise operator sensitivity
\begin{equation}\label{eq:Qop}
Q_K^{\mathrm{op}}=\tfrac12\sup_{\substack{\operatorname{supp}b\subset\operatorname{supp}\chi_K\\\norm{b}_{C^2(\Omega)}=1}}B_v(b,b).
\end{equation}
Then
\[
\tfrac12 B_v(e_K,e_K)\le Q_K^{\mathrm{op}}\norm{e_K}_{C^2}^2.
\]
$Q(K)$ of \Cref{eq:Q} samples this sensitivity in the probe direction $\beta_K$. 

Suppose a local resolution parameter $H_K$ produces the error estimate
\[
\norm{e_K}_{C^2}\le C_K H_K^r,\qquad r>0.
\]
Neglecting the cross-cell terms in \Cref{eq:regret-full} and replacing each diagonal contribution by its operator-norm bound leads to the allocation surrogate
\begin{equation}\label{eq:regret-additive}
\widetilde{\mathcal E}_{\mathrm{diag}}(H):=\sum_K Q_K^{\mathrm{op}} C_K^2 H_K^{2r}.
\end{equation}

\begin{proposition}\label{prop:alloc}
Let the computational work associated with the local resolutions be modeled by
\[
\sum_K H_K^{-d}\le B
\]
in $d$ spatial dimensions. Among cells satisfying $Q_K^{\mathrm{op}}C_K^2>0$, the resolution minimizing $\widetilde{\mathcal E}_{\mathrm{diag}}$ in \Cref{eq:regret-additive} under this budget satisfies
\begin{equation}\label{eq:alloc}
H_K\propto\bigl(Q_K^{\mathrm{op}}C_K^2\bigr)^{-1/(2r+d)}.
\end{equation}
Cells with zero weight do not require refinement under this surrogate.
\end{proposition}

\begin{proof}
At the optimum the work constraint is saturated. With multiplier $\Lambda>0$, stationarity of
\[
\sum_K Q_K^{\mathrm{op}}C_K^2H_K^{2r}+\Lambda\left(\sum_K H_K^{-d}-B\right)
\]
gives
\[
2r\,Q_K^{\mathrm{op}}C_K^2H_K^{2r-1}=\Lambda d\,H_K^{-d-1}.
\]
Hence
\[
H_K^{2r+d}=\frac{\Lambda d}{2r\,Q_K^{\mathrm{op}}C_K^2},
\]
which yields \Cref{eq:alloc}. Equivalently, with $z_K=H_K^{-d}$, the objective becomes a sum of terms proportional to $z_K^{-2r/d}$, which are convex for $z_K>0$; hence the stationary point is the global minimizer.
\end{proof}

\Cref{prop:alloc} is exact for the surrogate $\widetilde{\mathcal E}_{\mathrm{diag}}$, not for the full regret. Two approximations separate the two: the quadratic expansion itself is valid only in the regular perturbative regime, and the surrogate drops the cross-cell interactions in \Cref{eq:regret-full}. Thus \Cref{eq:alloc} should be interpreted as an asymptotic allocation principle rather than as a complete adaptive algorithm.

\subsection{Finite-amplitude and oracle diagnostics}
\label{subsec:oracle-diagnostics}

The quadratic sensitivity $Q(K)$ describes infinitesimal perturbations in the regular regime. To probe finite perturbations, we also use the \emph{finite-amplitude probe score}
\begin{equation}\label{eq:S}
S_{\mathrm{route}}(K)=\frac{\Reg_{\mathrm{route}}(v+\delta v_K;\,v)}{\norm{\delta v_K}},
\end{equation}
where $\delta v_K$ is localized near cell $K$, and $\norm{\delta v_K}=\bigl(|K|^{-1}\int_K|\delta v_K|^2\,\mathrm dx\bigr)^{1/2}$ is its root mean square over $K$. At a fixed, noninfinitesimal amplitude, $S_{\mathrm{route}}(K)$ is a diagnostic ranking score, not the infinitesimal $Q(K)$ of \Cref{eq:Q}. An executed-policy probe is a distinct quantity,
\begin{equation}\label{eq:S-pol}
S_{\mathrm{pol}}(K)=\frac{\Reg_{\mathrm{pol}}(v+\delta v_K;\,v)}{\norm{\delta v_K}}.
\end{equation}
The two scores quantify different notions of regret and need not agree.

For larger model changes, including changes that may cross a route-selection boundary, we use an oracle marginal quantity. Let $v_H$ denote a baseline reduced model and let $v_H^{(K)}$ denote the model obtained by refining cell $K$ while leaving the remaining cells unchanged. Define the \emph{oracle marginal value of resolution}
\begin{equation}\label{eq:marg}
\Delta_K^{\mathrm{marg}}=\Reg_{\mathrm{pol}}(v_H;v)-\Reg_{\mathrm{pol}}(v_H^{(K)};v).
\end{equation}
Positive values mean that refining cell $K$ reduces the executed regret; the quantity may be negative if refinement triggers a less favorable route selection. Unlike the quadratic approximation, \Cref{eq:marg} is defined directly through the nonlinear planning problem and therefore remains meaningful across branch switches. Its disadvantage is that it is fully oracle, requiring the true field $v$ and a refined-model evaluation.

The localized quantities introduced above play different roles and are not interchangeable: $Q(K)$ describes infinitesimal regular sensitivity, $S_{\mathrm{branch}}(K)$ the near-singular branch-ranking fragility, $S_{\mathrm{route}}(K)$ and $S_{\mathrm{pol}}(K)$ finite-amplitude regret-based scores, and $\Delta_K^{\mathrm{marg}}$ the nonlinear marginal value of refining a cell.

\begin{remark}\label{rem:oracle}
All of these are oracle diagnostics that require information about the true field $v$: $Q(K)$ uses the linearized route response and the second variation, whereas $S_{\mathrm{branch}}(K)$ uses the branch-action gap and the difference in directional branch-action derivatives. The finite-amplitude scores require perturbed-model planning and true-field regret evaluation, while $\Delta_K^{\mathrm{marg}}$ compares true-field regrets before and after cell refinement. A practical indicator should instead be computable from the reduced model alone---through a Jacobi/adjoint surrogate or a policy-margin quantity---and validated against these diagnostics; this is left for future work.
\end{remark}

\section{Numerical method and experiments}
\label{sec:numerics}

This section describes the numerical discretization and then illustrates the regular and singular stability mechanisms of \Cref{sec:dichotomy}, the localized sensitivity measures and resolution-allocation framework of \Cref{sec:sensitivity}, and the robustness of the singular branch switch under grid refinement.

\subsection{Numerical setup and Bellman discretization}
\label{sec:impl}

\paragraph{Bellman discretization}
All experiments use the same monotone semi-Lagrangian Bellman scheme. On the square domain $\Omega=[-1000\,\mathrm m,1000\,\mathrm m]^2$, the $n\times n$ uniform grid has nodes $x_i$ and spacing $h=2000/(n-1)\,$m. For a short look-ahead time $\tau$ and a heading set $\mathcal U_{N_A}$ of $N_A$ directions, the discrete value satisfies \begin{equation}\label{eq:sl}
V_h[w](x_i)=\min_{u\in\mathcal U_{N_A}}\left\{\tau+ \mathcal P_h[V_h[w]]\bigl(x_i+\tau(w(x_i)+s\,u)\bigr)\right\},\qquad x_i\notin\Gamma,
\end{equation} 
where $\mathcal P_h$ denotes bilinear interpolation~\cite{falcone2014semilagrangian,falcone2002semilagrangian}. Here $V_h[w]$ is the discrete value for planning field $w$, and $V_h[v]$ that of the true field. We solve \Cref{eq:sl} by fast-sweeping fixed-point iteration~\cite{zhao2005fast,sethian1996fast,sethian2003ordered}.

Write $\pi_{w,h}^*(x_i)\in\mathcal U_{N_A}$ for a minimizing heading selected in \Cref{eq:sl} at node $x_i$, choosing one if several minimize. For a grid policy $\pi$ with $\pi(x_i)\in\mathcal U_{N_A}$, define $J_h[w,\pi]$ as the grid arrival-time function satisfying the fixed-policy update
\[
J_h[w,\pi](x_i)=\tau+\mathcal P_h[J_h[w,\pi]]\bigl(x_i+\tau(w(x_i)+s\pi(x_i))\bigr),\qquad x_i\notin\Gamma,
\]
with $J_h[w,\pi](x_i)=0$ at target nodes and the same boundary treatment as in \Cref{eq:sl}. Consequently, value computation, policy extraction, and policy evaluation are discretely consistent at each grid node: 
\[
J_h[v,\pi_{v,h}^*](x_i)=V_h[v](x_i)
\]
to solver tolerance. Hence the discrete executed-policy regret
\[
\Reg_{\mathrm{pol},h}(\vhat;v,x_0)=J_h[v,\pi_{\vhat,h}^*](x_0)-V_h[v](x_0)\ge0
\] 
is nonnegative up to solver tolerance by construction. All reported discrete regrets and probe scores are evaluated at the specified initial point $x_0$; this argument is suppressed when fixed. This is an internal consistency property of the discrete operator, not an error bound for the continuous minimum-time problem: comparing a continuous trajectory arrival time with a discrete Bellman baseline can produce apparent negative regret through mismatched approximations.

\paragraph{Test fields and model reduction}
Both test fields are deterministic and available in closed form, so the reduction error is known exactly. The \emph{regular} shear field is decomposed as
\begin{equation}\label{eq:field-shear}
v=v^{\mathrm{base}}+v^{\mathrm{fine}},\qquad v^{\mathrm{base}}=\Bigl(U_{\mathrm{reg}}e^{-(y/\ell)^2},0\Bigr),\qquad v^{\mathrm{fine}}=\bigl(\partial_y\psi,-\partial_x\psi\bigr),
\end{equation}
where
\[
\psi=\frac{a}{k}\sin(k_xx+k_yy),\qquad (k_x,k_y)=k(\cos\theta,\sin\theta),\qquad k=\frac{2\pi}{\lambda}.
\]
The total field is rescaled, when necessary, so that $\abs{v}\le0.9\,s$.

The \emph{singular} field is the opposing westward jet
\[
v(x,y)=\Bigl(-U_{\mathrm{sing}}e^{-(y/\ell_j)^2},0\Bigr).
\]
Here $\ell_j$ is the opposing-jet width scale. To break its north--south symmetry we add a vertical bias to the $y$-component of the flow velocity 
\[
v_y\longmapsto v_y+\eta e^{-(y/200)^2}.
\]
Here $\eta$ is the signed vertical-bias amplitude in m/s, with an envelope of vertical scale $200\,\mathrm m$.

The reduced model is obtained by a componentwise Gaussian low-pass operator $R_\sigma$ whose kernel has standard deviation $\sigma$ grid cells (physical width $\sigma h$); values outside $\Omega$ are taken from the nearest boundary node. The resulting reduced field is
\[
v_H=R_\sigma v.
\]
The reduction error generally differs from the prescribed fine-scale component $v^{\mathrm{fine}}$. 

\Cref{tab:spec} summarizes the parameters used in the numerical experiments.

\begin{table}[t]
\centering
\caption{Parameters for the regular and singular fields, model reduction, discretization, and cellwise sensitivity and refinement experiments.}
\label{tab:spec}
\smallskip
\small
\setlength{\tabcolsep}{4pt}
\begin{tabular}{@{}lll@{}}
\toprule
group & quantity & value \\
\midrule
domain
& $\Omega$, own-speed $s$
& $[-1000\,\mathrm m,1000\,\mathrm m]^2$,\quad $1.5\,$m/s \\

& target set $\Gamma$
& $\overline B_{40}(x_T)$ \\

regular field
& jet $U_{\mathrm{reg}}$, width $\ell$
& $0.8\,$m/s,\quad $280\,$m \\

& fine-scale $a$, $\lambda$, $\theta$
& $0.25\,$m/s,\quad $180\,$m,\quad $35^\circ$ \\

& speed cap
& $0.9s=1.35\,$m/s \\

& start $x_0$, target center $x_T$
& $(-400,600)\,$m,\quad $(500,0)\,$m \\

singular field
& opposing-jet $U_{\mathrm{sing}}$, width $\ell_j$
& $1.3\,$m/s,\quad $110\,$m \\

& bias envelope, amplitude
& $e^{-(y/200)^2}$,\quad $\eta$ swept about $0$ \\

& start $x_0$, target center $x_T$
& $(-800,0)\,$m,\quad $(800,0)\,$m \\

reduction
& $R_\sigma$
& \begin{tabular}[c]{@{}l@{}}
Gaussian low-pass,\\
standard deviation $\sigma$ grid cells
\end{tabular}\\

& allocation baseline
& $\sigma=12$ grid cells \\

scheme
& grid $n\times n$
& $n\in\{121,161,201\}$ (by experiment) \\

& headings $N_A$
& $72$ \\

& look-ahead $\tau$, tolerance
& $1.5h/s$,\quad $10^{-6}$ \\

& robustness grids $n\times n$
& $n\in\{241,321,401\}$ \\

& robustness headings $N_A$
& $N_A\in\{48,96\}$ \\

& robustness look-ahead factors $\tau s/h$
& $\{1.0,2.0,3.0\}$ \\

cells
& refinement lattice
& $11\times11$ \\

& refinement-cell width
& $8$ grid cells per side \\

& probe $\delta v_K$
& \begin{tabular}[c]{@{}l@{}}
fixed-amplitude localized\\
velocity perturbation
\end{tabular} \\
\bottomrule
\end{tabular}
\end{table}

\paragraph{Implementation details and error measures}
On the uniform $n\times n$ grid we use
\[
\mathcal U_{N_A}=\left\{(\cos\theta_j,\sin\theta_j): \theta_j=\frac{2\pi j}{N_A},\quad j=0,\ldots,N_A-1\right\}
\]
and
\[
\tau=1.5\,h/s,
\]
so the self-propulsion contribution $s\tau$ corresponds to $1.5$ grid cells; the full Bellman displacement also depends on the flow. The target condition is $V_h[w](x_i)=0$ for $x_i\in\Gamma$. A foot point that falls outside $\Omega$ is projected onto $\Omega$ before interpolation. Rejecting such candidates instead leaves the true-field value and baseline regret of the regular experiment, and the tie and off-tie values and corridors of the singular experiment, unchanged to $10^{-6}$. Fast sweeping terminates when the maximum nodal change in one sweep is below $10^{-6}$, after at least three complete sweeps.

The Gaussian reduction is truncated at $4\sigma$. The speed cap in \Cref{eq:field-shear} is applied before reduction through the global rescaling
\[
v\longmapsto v\,\min\left\{1,\frac{0.9s}{\max_\Omega\abs{v}}\right\}.
\]

All discrete field and value errors use the root-mean-square norm
\begin{equation}\label{eq:discrete-rms}
\norm{f}_{2,h;\mathcal N}=\left(\frac{1}{|\mathcal N|}\sum_{x_i\in\mathcal N}\abs{f(x_i)}^2\right)^{1/2}
\end{equation}
over the relevant node set $\mathcal N$. For a refinement cell $K$ we abbreviate this by $\norm{\cdot}_{2,h;K}$.

For the scores of \Cref{eq:S,eq:S-pol}, the localized probe is
\[
\delta v_K=g_K\,(\partial_y\psi_{\mathrm p},-\partial_x\psi_{\mathrm p}), 
\]
where
\[
\psi_{\mathrm p}=\frac{A_{0}}{k_{\mathrm{p}}}\sin(k_{x,\mathrm{p}}x+k_{y,\mathrm{p}}y),\qquad g_{K}(x)=\exp\left(-\frac{|x-x_{K}|^{2}}{2r_{K}^{2}}\right), 
\]
with $x_K$ the center of cell $K$ and $r_K=1.3\times 8h=10.4\,h$ (about $173\,$m on the $121\times121$ grid). We fix $A_0=0.12\,$m/s, wavelength $2\pi/k_{\mathrm p}=150\,$m, and wave-vector angle $40^\circ$ across all cells, with $(k_{x,\mathrm p},k_{y,\mathrm p})=k_{\mathrm p}(\cos40^\circ,\sin40^\circ)$; the stream-function derivatives are taken on the grid, and every probed field is verified strictly subcritical. The norm $\norm{\delta v_K}$ in \Cref{eq:S} is realized numerically by $\norm{\delta v_K}_{2,h;K}$.

\paragraph{Route and policy probe evaluation}
The retained scores use the consistent Bellman operator: $S_{\mathrm{pol},h}(K)$ denotes the executed-policy score
\begin{equation}\label{eq:policy-score-discrete}
S_{\mathrm{pol},h}(K)=\frac{J_h[v,\pi_{v+\delta v_K,h}^*](x_0)-V_h[v](x_0)}{\norm{\delta v_K}_{2,h;K}},
\end{equation}
which is the quantity used in the sensitivity maps and correlations.

The $11\times11$ refinement cells have centers on a regular subgrid, each cell spanning $8\times8$ grid cells, and are disjoint. The controlled rate test in \Cref{subsec:num-regular} uses a $201\times201$ grid, the sensitivity and allocation experiments in \Cref{subsec:num-allocation} use a $121\times121$ grid, and the singular-regime experiments in \Cref{sec:sweep} use a $161\times161$ grid. The default heading count is $N_A=72$.

The code and data that generate all figures and tables are available from the corresponding author upon reasonable request.

\subsection{Regular regime: regret scaling}
\label{subsec:num-regular}

\Cref{thm:regular} predicts that, for a sufficiently small smooth model perturbation, the route displacement is first order while the committed-route and executed-policy regrets are second order. 

\begin{figure}[t]
\centering
\includegraphics[width=\textwidth]{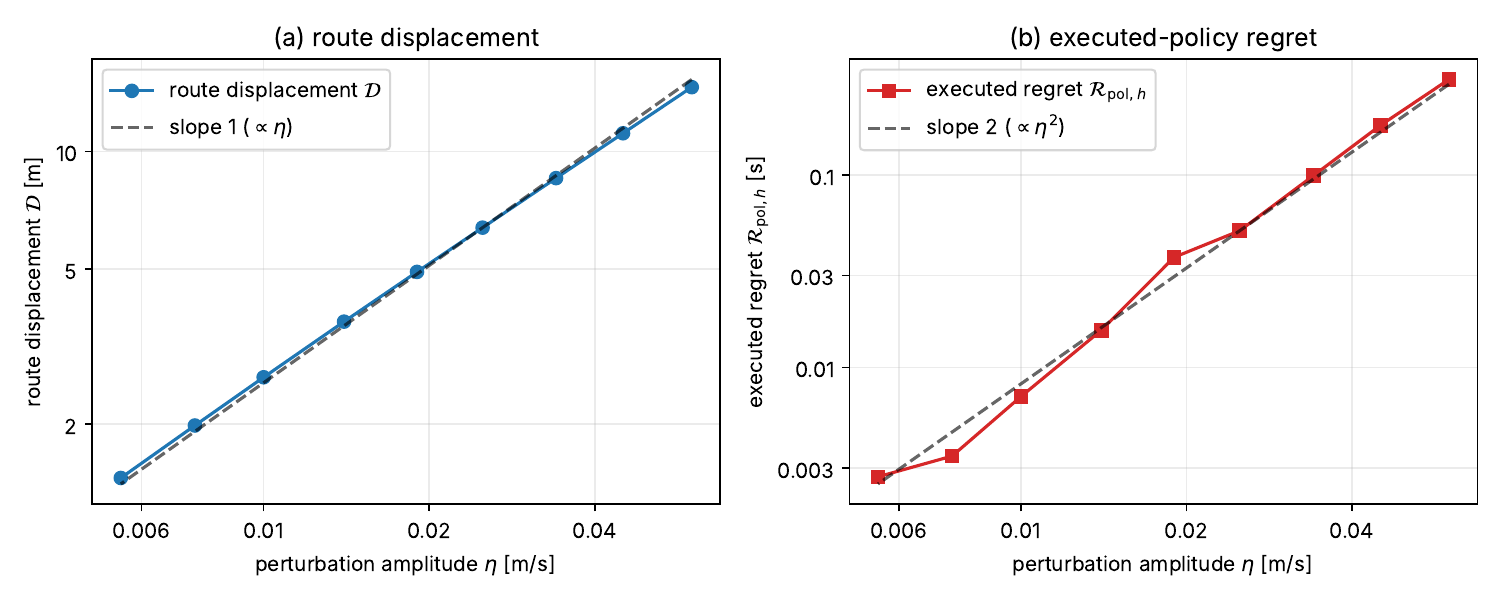}
\caption{Controlled single-direction rate test, $w_\eta=v+\eta\beta$ for a fixed smooth $\beta$. (a) Route displacement $\Sel$ follows the slope-one guide ($\Sel\propto\eta$), consistent with the first-order route-sensitivity bound. (b) Executed-policy regret $\Reg_{\mathrm{pol},h}$ follows the slope-two guide ($\Reg_{\mathrm{pol},h}\propto\eta^2$). Fitted exponents are $0.96$ and $2.08$.}
\label{fig:rate}
\end{figure}

\paragraph{Controlled rate test}
To test the predicted small-perturbation scaling, we fix a single smooth divergence-free direction $\beta=(\partial_y\psi_\beta,-\partial_x\psi_\beta)$ with $\psi_\beta\propto\sin\!\big(k(x\cos55^\circ$ $+y\sin55^\circ)\big)$, $k=2\pi/700\,\mathrm m^{-1}$, normalized to $\max_\Omega\abs{\beta}=1$, set $w_\eta=v+\eta\beta$, and sweep the scalar amplitude $\eta\in[0.0055,0.06]\,$m/s. For each $\eta$ we record the route displacement $\Sel(w_\eta;v)$, measured between the route committed in $w_\eta$ and the true optimal route, and the executed-policy regret $\Reg_{\mathrm{pol},h}(w_\eta;v)$. The two quantities use different constructions. Routes are reconstructed by integrating the closed-loop dynamics with heading $-\nabla V_h/\abs{\nabla V_h}$, using the bilinearly interpolated gradient of the discrete value and a time step of $0.25\,$s; this continuous reconstruction allows route differences below the grid spacing ($h=10\,$m here) to be measured, although such measurements remain subject to the spatial discretization error. The regret $\Reg_{\mathrm{pol},h}$ is instead evaluated with the Bellman-optimal policy on the finite heading set, as in \Cref{sec:impl}. \Cref{fig:rate} plots both against $\eta$ on logarithmic axes. A least-squares fit over the tested amplitudes gives log--log slopes of $0.96$ for $\Sel$ and $2.08$ for $\Reg_{\mathrm{pol},h}$, in agreement with the first-order route displacement of \Cref{thm:regular} and the second-order executed-policy regret of \Cref{lem:track}. Across other grids and heading sets (\Cref{tab:rate-slopes}), the displacement slope stays at $0.95$--$0.97$, while the regret slope ranges from $1.79$ to $2.23$ depending on the discretization and the fitting range; the regrets in this test are $0.002$--$0.3\,$s, close to the resolution of the discrete policy.

\begin{table}[t]
\centering
\caption{Log--log fitted slopes for route displacement $\Sel$ and executed-policy regret $\Reg_{\mathrm{pol},h}$ across grid sizes and heading counts. Each entry reports the displacement/regret slopes. The reference setting is $n=201$ and $N_A=72$.}
\label{tab:rate-slopes}
\smallskip
\small
\begin{tabular}{@{}cccc@{}}
\toprule
grid $n$ & headings $N_A$ & all amplitudes & excluding two smallest amplitudes \\
\midrule
$201$ & $72$ & $0.96/2.08$ & $0.96/2.08$ \\
$161$ & $72$ & $0.97/2.09$ & $0.96/2.07$ \\
$241$ & $72$ & $0.96/1.95$ & $0.95/1.79$ \\
$201$ & $48$ & $0.97/2.15$ & $0.96/2.16$ \\
$201$ & $96$ & $0.96/2.23$ & $0.95/2.13$ \\
\bottomrule
\end{tabular}
\end{table}

\subsection{Localized sensitivity and resolution allocation}
\label{subsec:num-allocation}

We next test whether local flow-error magnitude identifies the cells that matter most to the navigation decision. Let $v_H$ be the baseline reduced field with $\sigma=12$ grid cells and define
\begin{equation}\label{eq:EH}
E_H(K)=\norm{v-v_H}_{2,h;K}.
\end{equation}
The probe score $S_{\mathrm{pol},h}(K)$ uses the same amplitude in every cell.

\begin{figure}[t]
\centering
\includegraphics[width=\textwidth]{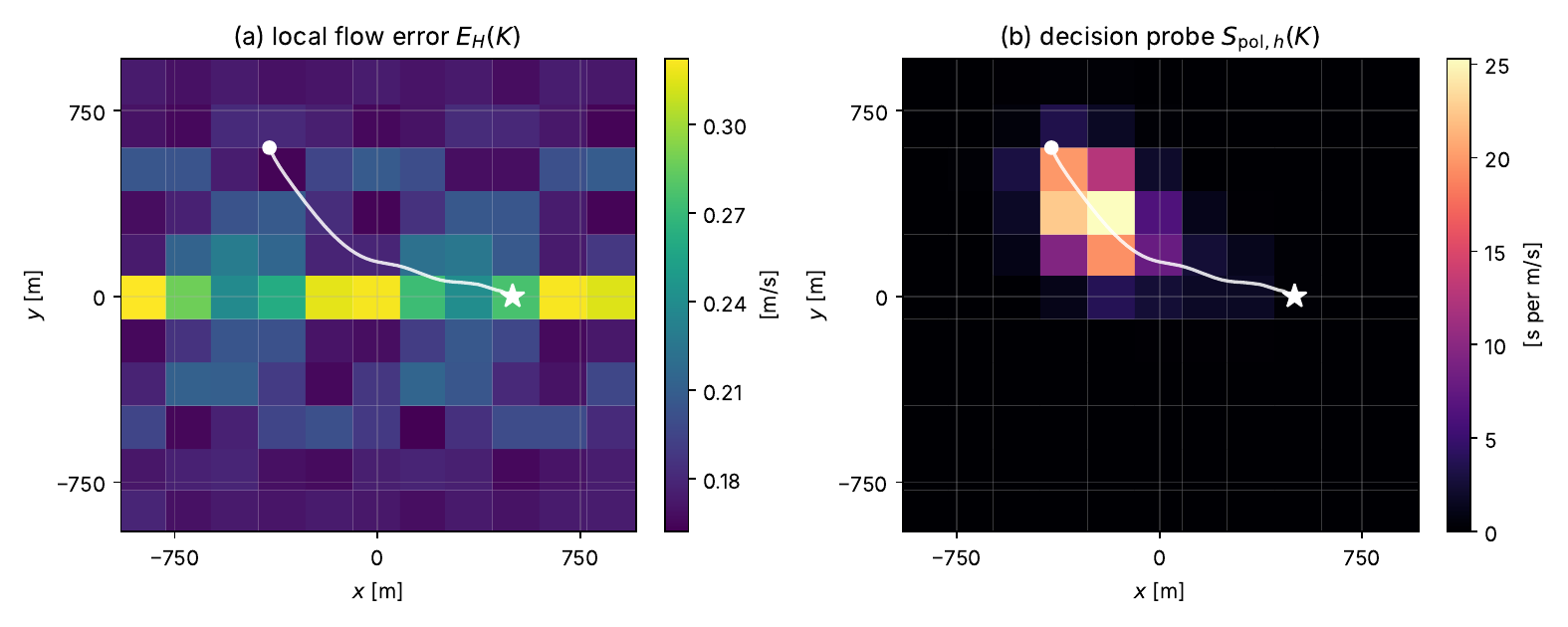}
\caption{Regular shear experiment. (a) Local flow-error magnitude $E_H(K)=\norm{v-v_H}_{2,h;K}$. (b) Finite-amplitude executed-policy probe score $S_{\mathrm{pol},h}(K)$. The true optimal route is overlaid in white. Large local flow error does not by itself identify the cells with the largest decision response.} 
\label{fig:sens}
\end{figure}

\paragraph{Flow error versus decision sensitivity}
\Cref{fig:sens} compares the spatial maps of $E_H(K)$ and $S_{\mathrm{pol},h}(K)$. The local flow-error magnitude is heterogeneous across the domain, but the cells with the largest flow errors do not coincide systematically with those producing the largest finite-amplitude decision-regret response. The true optimal route is overlaid to indicate the geometry of the decision-sensitive region.

Compared cell by cell over the $121$ refinement cells, the two maps have a Pearson correlation of $0.02$ and a Spearman rank correlation of $0.35$. These values suggest that, in this experiment, local flow-error magnitude is a poor proxy for the decision score: the largest flow errors lie along the jet axis, whereas the decision-sensitive cells cluster near the start of the committed route.

\paragraph{Oracle allocation benchmark}
The asymptotic allocation rule of \Cref{prop:alloc} applies in the regular perturbative regime and depends on the quadratic sensitivity. The finite cell refinements used here can be much larger and can cross route-selection boundaries. In the cellwise allocation experiments, refining a cell means replacing $v_H$ by $v$ on that cell, leaving the others unchanged; this hard replacement produces jumps at cell interfaces and serves as a discrete oracle benchmark, to which the $C^2$ bounds of \Cref{thm:regular} do not directly apply. A continuum comparison with that theorem requires smooth blending, carried out at the streamfunction level if incompressibility is needed. To rank cells in this benchmark, define the marginal-regret reduction
\[
\Delta_{K,h}^{\mathrm{marg}}:=\Reg_{\mathrm{pol},h}(v_H;v,x_0)-\Reg_{\mathrm{pol},h}(v_H^{(K)};v,x_0).
\]
This is the discrete counterpart of $\Delta_K^{\mathrm{marg}}$ in \Cref{eq:marg}. 

We therefore benchmark the \emph{potential} benefit of decision-aware allocation using the nonlinear oracle quantity $\Delta_{K,h}^{\mathrm{marg}}$.

For each baseline cell, $\Delta_{K,h}^{\mathrm{marg}}$ is computed once by applying this oracle refinement to that cell and reevaluating the executed-policy regret. The ranking is not recomputed after subsequent refinements. Thus, the benchmark requires one additional planning/evaluation calculation per candidate cell, an $O(N)$ oracle family for $N$ cells, and intentionally ignores interactions that would be captured by re-ranking after each refinement.

At equal refinement budgets, we compare three strategies: \emph{nested uniform} refinement, \emph{largest-flow-error-first refinement}, which ranks cells by $E_H(K)$, and \emph{oracle marginal-regret refinement}, which ranks them by $\Delta_{K,h}^{\mathrm{marg}}$. Nested uniform refinement uses the farthest-point ordering of the cells, so each budget's cells are contained in the next; ties are broken by cell index. The budgets are $\{0,4,8,16,28,44,64\}$ refined cells.

\begin{figure}[t]
\centering
\includegraphics[width=0.66\textwidth]{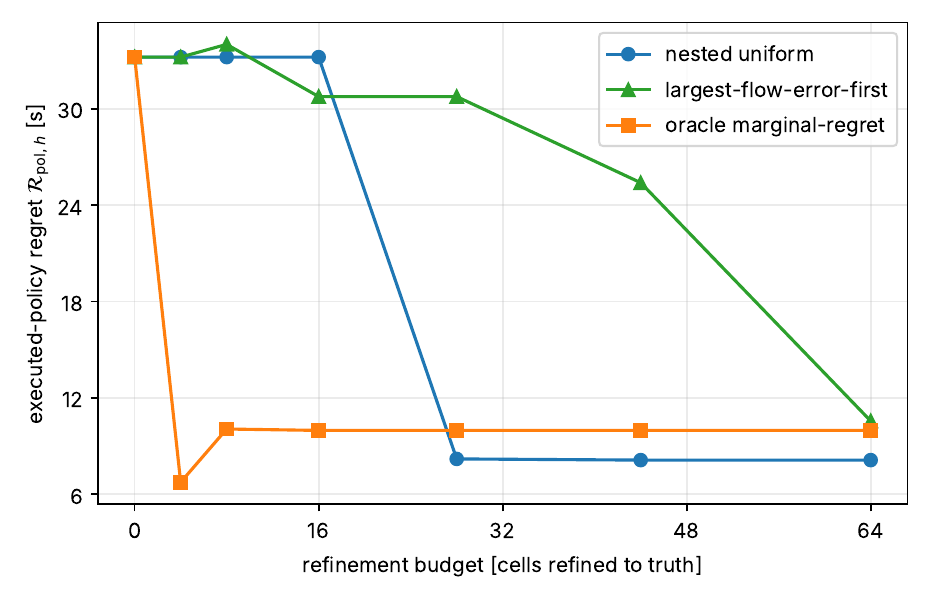}
\caption{Executed-policy regret versus refinement budget for three cell-allocation strategies. The oracle marginal-regret ranking reaches low regret with only four refined cells, while largest-flow-error-first refinement and nested uniform refinement need many more. The oracle curve rises from its budget-four minimum to a plateau under this precomputed oracle ranking, which is not re-evaluated as cells are added.}
\label{fig:adapt}
\end{figure}

\Cref{fig:adapt} shows that the oracle marginal-regret ranking reduces the executed regret from approximately $33\,$s to about $7\,$s after only four refined cells. Largest-flow-error-first refinement leaves the regret at $31$--$34\,$s through budget $28$, and nested uniform refinement remains at the baseline value through budget $16$; both require roughly $28$--$64$ cells to reach a comparable regret. The advantage of the oracle ranking is therefore its efficiency at small budgets, not a uniform margin at every budget.

The finite-amplitude probe $S_{\mathrm{pol},h}(K)$ is only a partial surrogate for the oracle marginal value. Across the $121$ cells, the rankings of $S_{\mathrm{pol},h}(K)$ and $\Delta_{K,h}^{\mathrm{marg}}$ have Spearman correlation $0.30$ and share about half of their top eight cells. Only $5$ of the $121$ cells (about $4\%$) have marginal value above $5\%$ of the maximum. These observations reinforce the distinction made in \Cref{subsec:oracle-diagnostics}: $S_{\mathrm{pol},h}(K)$ is a diagnostic finite-amplitude probe, whereas $\Delta_{K,h}^{\mathrm{marg}}$ measures the nonlinear marginal value of a specific refinement.

To check that the comparison in \Cref{fig:adapt} is not specific to the $121\times121$ grid, we repeated it on a $241\times241$ grid with the physical filter width ($200\,$m), the refinement cells ($133\,$m, same centers), and the budgets unchanged. The baseline regret is $35.5\,$s; the oracle ranking reduces it to $7.3\,$s at budget $4$ and to $10.4\,$s at larger budgets, largest-flow-error-first refinement stays at $32$--$36\,$s through budget $44$, and nested uniform refinement stays at the baseline through budget $16$ and reaches $8.1\,$s at budget $28$. The oracle marginal values on the two grids have Spearman correlation $0.93$, share the same top four cells, and again exceed $5\%$ of the maximum in only $5$ cells.

\paragraph{Refinement and monotonicity}
\label{sec:branchroutes}
In this experiment, the nested uniform refinement sequence yields a monotone decrease in executed regret, from $33\,$s at the baseline to about $8\,$s by budget $28$. Nestedness guarantees that the discrete field error is nonincreasing, but it does \emph{not} guarantee monotone regret: reducing the field error can disturb cancellations between errors or change the selected branch. Indeed, the oracle sequence in \Cref{fig:adapt}, whose top-$k$ cells are also nested, is itself non-monotone, rising from about $7\,$s at budget $4$ to about $10\,$s at budget $8$.

\subsection{Singular regime: branch switching and stability radius}
\label{sec:sweep}

We now illustrate \Cref{prop:switch,cor:radius}. The opposing-jet field has two symmetric north/south routes at zero bias. A small bias breaks the tie.

\begin{figure}[t]
\centering
\includegraphics[width=\textwidth]{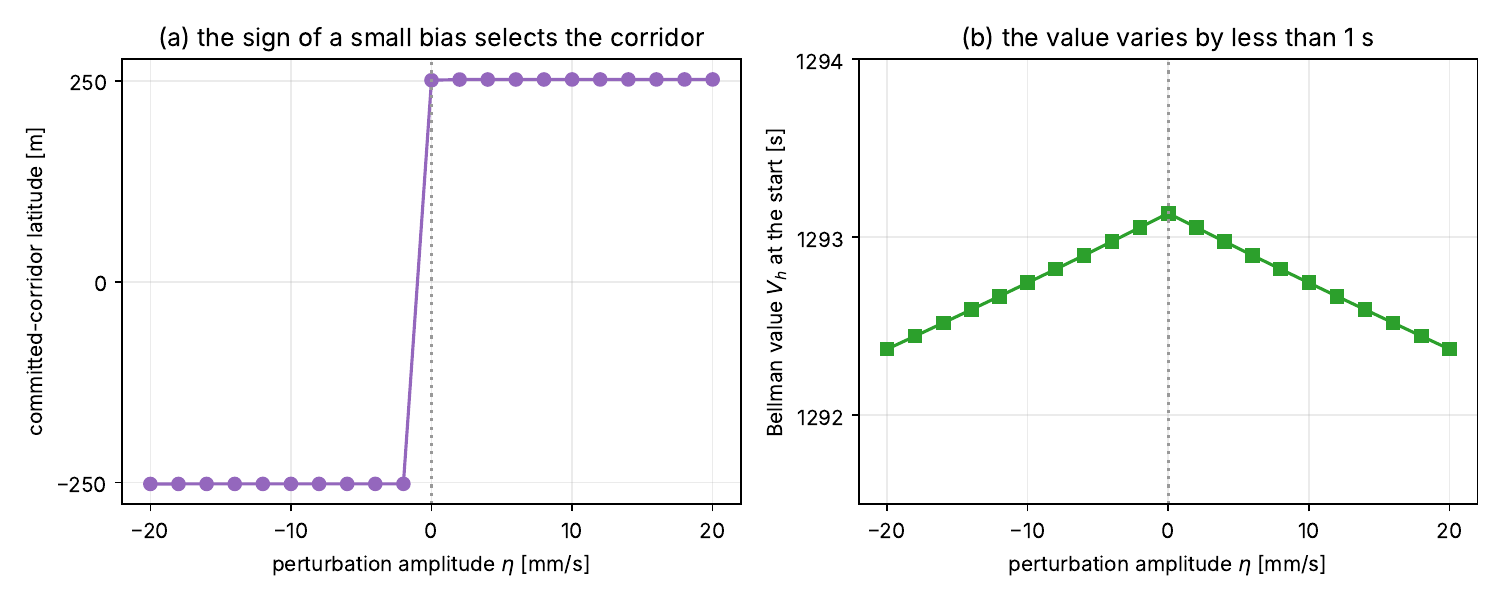}
\caption{Singular opposing-jet experiment. (a) The sign of a tiny asymmetric perturbation selects the north or south corridor, and the committed route changes by about $503\,$m across zero perturbation. (b) The minimum-time estimate, evaluated with the consistent Bellman operator, is $1293.13\,$s at the tie and decreases to $1292.37\,$s at $|\eta|=20\,$mm/s, a variation of $0.76\,$s ($0.06\%$).}
\label{fig:switch}
\end{figure}

\paragraph{Value stability and route discontinuity}
\Cref{fig:switch} shows the central singular effect. Opposite signs at the smallest nonzero amplitude, $\abs{\eta}=2\,$mm/s, already select opposite corridors. The selected route switches from the southern corridor ($y\approx-251.7\,$m) to the northern one ($y\approx+251.7\,$m), a lateral change of about $503\,$m.

The Bellman value remains near $1293\,$s throughout the bias sweep. More importantly for the branch-switch claim, it changes by less than $0.1\,$s between the tie and the adjacent sweep amplitudes $\eta=\pm2\,$mm/s, while the selected route changes by roughly half a kilometer. Thus the discontinuity appears in the selected minimizer rather than in the minimum value, in agreement with \Cref{prop:switch}.

\paragraph{Exact and near-degenerate branch switching}
To resolve the switching structure more directly, we sweep the bias over $21$ equally spaced values. For the symmetric problem, $\eta\in[-20,20]\,$mm/s (\Cref{fig:switch}); the committed corridor stays near $-251.7\,$m for $\eta<0$ and near $+251.7\,$m for $\eta>0$, without an intermediate branch.

\begin{figure}[t]
\centering
\includegraphics[width=0.58\textwidth]{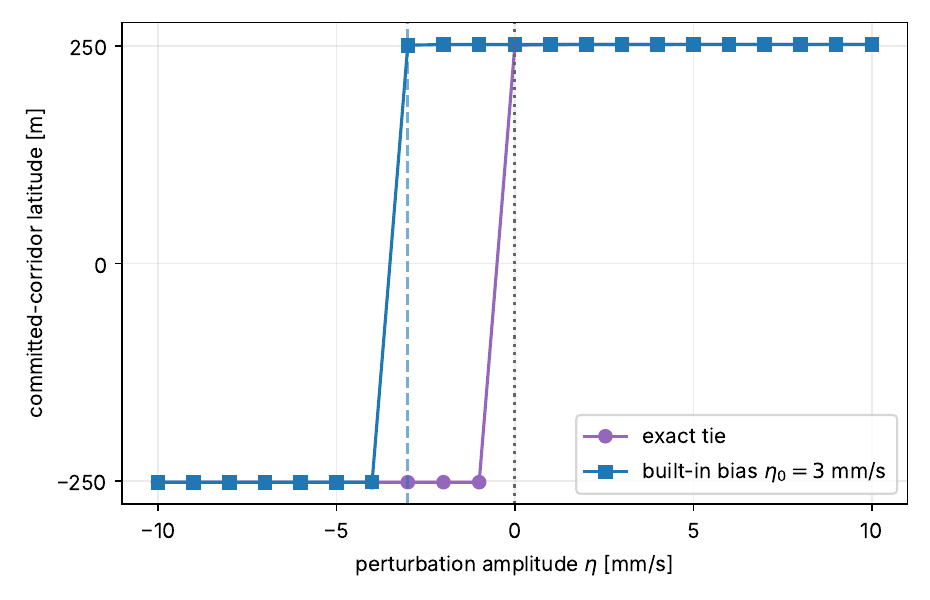}
\caption{Committed-corridor latitude versus perturbation amplitude $\eta$ in the singular opposing-jet experiment. Exact tie (circles): the selected corridor stays near its respective branch on either side and jumps by about $503\,$m through $\eta=0$, consistent with \Cref{prop:switch}. Near-degenerate case with built-in velocity bias $\eta_0=3\,$mm/s (squares): the step moves to $\eta\approx-\eta_0$, where the same-direction biases cancel. This illustrates a finite stability interval; the coefficient in \Cref{eq:radius} is tested independently in \Cref{fig:radius}. The parameter $\eta_0$ is a velocity bias, not the branch-action gap. Both curves use $21$ equally spaced amplitudes in $[-10,10]\,$mm/s.}
\label{fig:ampsweep}
\end{figure}

We then add a permanent bias $\eta_0=3\,$mm/s, which to leading order creates a branch-action gap $\Delta=\eta_0\abs{G_\beta(\gamma_1)-G_\beta(\gamma_2)}+O(\eta_0^2)$. Sweeping an additional perturbation over $\eta\in[-10,10]\,$mm/s then shifts the switch to $\eta\approx-\eta_0$ (\Cref{fig:ampsweep}): with both biases in the same direction $\beta$ the total field is $v+(\eta_0+\eta)\beta$, so symmetry is restored at $\eta=-\eta_0$ by construction.

\begin{figure}[t]
\centering
\includegraphics[width=\textwidth]{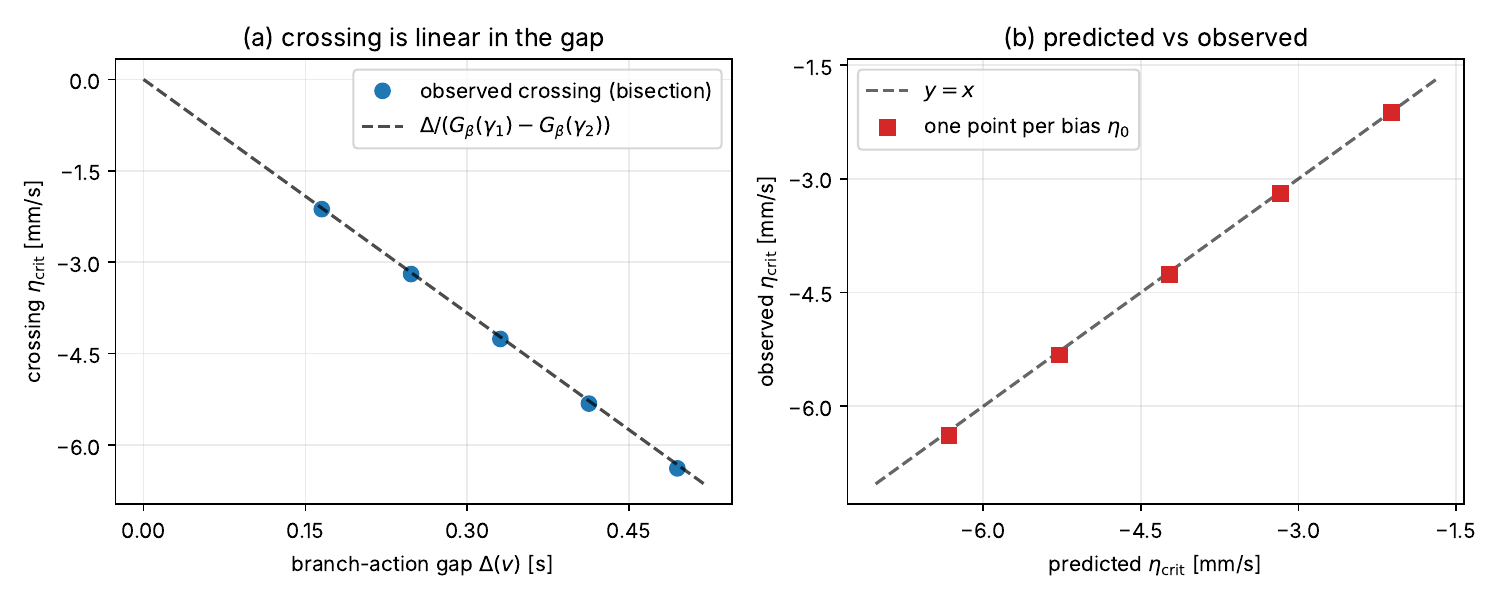}
\caption{Independent test of the stability-radius coefficient \Cref{eq:radius}. The branch-action gap $\Delta$ and the slope difference $G_\beta(\gamma_1)-G_\beta(\gamma_2)$ are measured from half-plane branch actions, independently of the crossing. The observed crossing is the amplitude at which the committed route of the full-domain solve changes corridor, located by bisection. (a) The observed crossing amplitude $\eta_{\mathrm{crit}}$ is linear in $\Delta$ and follows the predicted slope $1/(G_\beta(\gamma_1)-G_\beta(\gamma_2))$. (b) Predicted versus observed $\eta_{\mathrm{crit}}$ lie near the diagonal, agreeing to within $1\%$.}
\label{fig:radius}
\end{figure}

\paragraph{Independent test of the radius coefficient}
To test the leading-order coefficient in \Cref{eq:radius}, we isolate the two branch actions by solving the minimum-time problem separately on each half-plane ($y\ge0$ for the north branch, $y\le0$ for the south branch), which yields $\Jf_v(\gamma_1)$ and $\Jf_v(\gamma_2)$ directly. On the symmetric field, the confined routes reach an apex of $251.8\,$m and stay at least $120\,$m from the artificial boundary $y=0$ except near their endpoints, so that constraint is inactive in the interior and the half-plane value is the genuine branch action. The bias and probe use distinct vertical-velocity profiles $\beta_0:\delta v_y=e^{-(y/200)^2}$ and $\beta:\delta v_y=e^{-(y/320)^2}$. For a near-degenerate family built from $\eta_0\beta_0$ with $\eta_0\in\{2,\ldots,6\}\,$mm/s, we measure the gap $\Delta=\Jf_v(\gamma_2)-\Jf_v(\gamma_1)$ from the half-plane values and the slope difference $G_\beta(\gamma_1)-G_\beta(\gamma_2)$ by a central difference of step $\pm2\,$mm/s in the direction $\beta$. The leading-order prediction $\eta_{\mathrm{crit}}=\Delta/(G_\beta(\gamma_1)-G_\beta(\gamma_2))$ is compared with two observations on the $161\times161$ grid, each located by bisection to $0.005\,$mm/s: the zero of the half-plane action difference, and the amplitude at which the committed route of the full-domain solve changes corridor. The first differs from the prediction by at most $0.009\,$mm/s ($0.15\%$), comparable to the bisection tolerance; since the prediction is built from the same half-plane values, this checks only that the action difference is linear in $\eta$ over the tested range. The second is independent of the half-plane solves and agrees to within $1\%$ (\Cref{fig:radius}). The two observed thresholds differ by $0.015$--$0.05\,$mm/s, which exceeds the bisection tolerance and reflects the difference between the half-plane and full-domain discrete computations. The test thus supports the predicted leading-order coefficient but does not resolve the order of the remainder.

\subsection{Discretization robustness}
\label{subsec:num-convergence}

Finally, we check that the singular observation is not an artifact of the $161\times161$ grid. \Cref{tab:conv} recomputes the value at the start on three grids, at the exact tie and at a small positive bias.

The absolute Bellman value changes by approximately $1$--$3\,$s over these meshes. The quantity relevant to the singular claim is much more stable: the difference $\abs{V_h[v_{\eta_b}](x_0)-V_h[v](x_0)}$ between the value at the tie and the value at bias $\eta_b=4\,$mm/s (the same for $-\eta_b$ by symmetry) is $0.15$, $0.16$, and $0.10\,$s on the three grids, whereas the committed-corridor separation remains between $502$ and $504\,$m. On the $161\times161$ grid, changing the heading set from $72$ to $48$ or $96$ directions changes the value at the tie by less than $0.6\,$s, and leaves this value difference ($0.15$--$0.16\,$s) and the corridor separation ($502$--$504\,$m) essentially unchanged. These mesh-to-mesh changes do not bound the absolute error relative to the continuous problem.

\begin{table}[t]
\centering
\caption{Grid refinement for the singular opposing-jet experiment. Here $v_{\eta_b}$ denotes the opposing-jet field with the vertical bias amplitude $\eta_b=+0.004\,\mathrm{m/s}$, whereas $v$ denotes the zero-bias field. The values at the exact tie and with this bias differ by only $0.10$--$0.16\,$s on each grid, while the committed-corridor separation remains $502$--$504\,$m. The absolute Bellman value varies by approximately $1$--$3\,$s across the grids.}
\label{tab:conv}
\smallskip
\begin{tabular}{ccccc}
\toprule
grid $n\times n$
& $V_h[v]$ at tie [s]
& $V_h[v_{\eta_b}]$ off tie [s]
& difference [s]
& corridor jump [m] \\
\midrule
$121\times121$ & $1296.336$ & $1296.184$ & $0.152$ & $504$ \\
$161\times161$ & $1293.135$ & $1292.976$ & $0.159$ & $503$ \\
$201\times201$ & $1294.172$ & $1294.076$ & $0.096$ & $502$ \\
\bottomrule
\end{tabular}
\end{table}

An independent continuous feasible-route calculation for the same zero-bias parameters---integrating the exact characteristic ODE system
\[
\dot x=f(y)+s\cos\theta,\qquad \dot y=s\sin\theta,\qquad \dot\theta=-f'(y)\cos^2\theta,\qquad f(y)=-U_{\mathrm{sing}}e^{-(y/\ell_j)^2},
\]
from $(-800,0)$ over the fan of extremals that return to the jet, and taking the smallest first-arrival time at the target disc---gives $1283.6185\,$s. This is a feasible upper bound for the continuous optimal time, not a proof of global optimality. The $161\times161$-grid Bellman value $1293.1\,$s exceeds that feasible arrival time by about $9.5\,$s, so the $1$--$3\,$s variation between these meshes should not be interpreted as an absolute error estimate. The excess is not monotone over the three grids of \Cref{tab:conv} ($12.7$, $9.5$, and $10.6\,$s for $n=121$, $161$, and $201$), but it decreases on finer grids, to $8.8\,$s ($241\times241$), $6.5\,$s ($321\times321$), and $5.0\,$s ($401\times401$), at approximately first order in $h$ for $n\ge201$. Changing the look-ahead factor in $\tau=1.5\,h/s$ to $1.0$, $2.0$, or $3.0$ on the $161\times161$ grid changes the value at the tie by at most $3.7\,$s, while the value difference between the tie and bias $\eta_b$ stays at $0.14$--$0.17\,$s and the corridor separation at $503$--$505\,$m.

These refinement tests support the qualitative robustness of the singular phenomenon: the minimum value at the tie and at bias $\eta_b$ differ by at most $0.17\,$s in every setting tested, even though the selected route changes by roughly half a kilometer.

\section{Conclusion}
\label{sec:disc}

This work studies how model reduction affects time-optimal navigation in multiscale flows. The main conclusion is that value accuracy alone does not determine the stability of the navigation decision. When the optimal route is unique and nondegenerate, small flow perturbations produce a smooth, first-order change in the route while the performance loss is second order; near a competition between distinct routes, the travel time can stay stable even as the selected route changes abruptly, a small branch-action gap providing only a finite stability margin. This also changes how model resolution should be allocated: in the regular regime, local errors can be weighted by their effect on the decision, giving an asymptotic resolution-allocation principle, whereas near competing routes, sensitivity is governed by the branch-action gap and by how perturbations change the relative branch values. The experiments show that local flow-error magnitude can be a poor guide to decision-relevant resolution, while an oracle marginal-regret refinement strategy achieves substantially lower regret. The results do not yet provide a practical adaptive algorithm; developing computable decision indicators and testing the framework in more realistic and multibranch settings are natural next directions.

\appendix
\section{Uniform feedback stability on a fixed tube}
\label{app:feedback}

We prove the feedback estimate \Cref{eq:feedback-stability} of \Cref{lem:track}: persistence of a single route from \Cref{lem:branch} does not control $V_w$ or its feedback on a neighborhood, so we build both on a fixed exterior tube $U$, uniformly for $w$ near $v$.

\emph{Characteristic coordinates.} Let $z(\sigma)$ parameterize a short arc of $\partial\Gamma$ with outward normal $n(z)$, and for $w$ near $v$ set $q_w(z)=n(z)/(s-w(z)\cdot n(z))$, so $\Ham_w(z,q_w)=1$. Integrating
\[
\dot Z=s\,P/\abs{P}-w(Z),\quad \dot P=Dw(Z)^\top P,\quad (Z,P)|_{t=0}=(z(\sigma),q_w(z(\sigma)))
\]
in backward target time $t$ defines $E_w(\sigma,t)=Z_w(\sigma,t)$.

\emph{Uniform invertibility.} At $w=v$, (A2) makes $E_v$ injective with nonsingular derivative on a compact strip around the reference route, where $P_v=\grad V_v(E_v)$ and $V_v(E_v)=t$. Since the characteristic flow is $C^1$ in $(w,\sigma,t)$ for the $C^2$ topology, the parametric inverse-function theorem inverts $E_w$ on a fixed smaller tube $U$, uniformly for $w$ near $v$; its time coordinate $W_w$ solves $\Ham_w(x,\grad W_w)=1$ with $W_w|_{\partial\Gamma}=0$.

\emph{Global identification.} $W_w$ is the global value on $U$: if not, take $w_j\to v$ and points $x_j$ in a compact substrip with $W_{w_j}(x_j)\ne V_{w_j}(x_j)$; the global minimizers are uniformly Lipschitz (constant-speed gauge, strict subcriticality), and a uniform limit is a true reference minimizer (lower semicontinuity of the convex length) whose terminal data lie in the strip, so $E_{w_j}$-invertibility identifies it with the constructed route---a contradiction, a length bound handling the target collar. Hence $W_w=V_w$ on $U$.

\emph{Uniform bounds.} The inverses and covectors depend $C^1$ on $w$, giving $$\norm{\grad V_{\vhat}-\grad V_v}_{C^0(U)}\le C\norm{\vhat-v}_{C^2}.$$ With $(s+\norm{w}_{L^\infty})^{-1}\le\abs{\grad V_w}\le c_0^{-1}$, normalizing yields the heading bound, which is \Cref{eq:feedback-stability}. On a slightly larger tube, $v\cdot n-s\le-c_0$ keeps the reference trajectory strictly inside $\Gamma$ at $V_v(x_0)+\delta$, so the lateral margin excludes exit and forces arrival---the uniform confinement and arrival-time bound used in \Cref{lem:track}.

\bibliographystyle{siamplain}
\bibliography{references}

\end{document}